\documentclass[12pt]{amsart}
\usepackage{amssymb,amsthm,amsmath}

\begin{document}
\theoremstyle{plain}
\newtheorem{thm}{Theorem}[section]
\newtheorem{lem}[thm]{Lemma}
\newtheorem{prop}[thm]{Proposition}
\newtheorem{cor}[thm]{Corollary}
\theoremstyle{definition}
\newtheorem{Def}[thm]{Definition}
\newtheorem{rem}[thm]{Remark}
\newtheorem{conj}[thm]{Conjecture}
\newtheorem{cl}{Claim}
\newtheorem{ex}[thm]{Example}
\def\Pf{\trivlist\item[\hskip\labelsep\textit{Proof.}]}
\def\endPf{\strut\hfill\framebox(6,6){}\endtrivlist}

\title[Bogomolov's conjecture for hyperelliptic curves]{Bogomolov's conjecture 
for hyperelliptic curves \\ over function fields}
\author{Kazuhiko Yamaki}
\date{January 27, 1999}
\address{Department of Mathematics, Faculty of Science,
Kyoto University, Kyoto,Japan}
\email{yamaki@kusm.kyoto-u.as.jp}
\maketitle

\section{Introduction}

Let us fix a field $k$.
Let $X$ be a smooth projective surface over $k$, 
$Y$ a smooth projective curve over $k$, 
and let $f : X \rightarrow Y$ 
be a generically smooth semistable curve of genus $g \geq 2$ over $Y$.
Let $K$ be the function field of $Y$, 
$\overline{K}$ the algebraic closure of $K$, 
and let $C$ be the generic fiber of $f$.
For $D \in \mbox{Pic}^{1}(C)(\overline{K})$, let
\[ j : C(\overline{K}) \rightarrow \mathrm{Pic}^{1}(C)(\overline{K}) \]
be a morphism defined by $j(x) = x - D$, and $\| \cdot \|_{NT}$
the semi-norm arising from the N\'{e}ron-Tate pairing on 
Pic$^{1}(C)(\overline{K})$. We set
\[ B_{C}(P;r) = \{ x \in C(\overline{K}) \mid \| j(x) - P \|_{NT} \leq r \} \]
for $P \in$ Pic$^0(C)(\overline{K})$ and $r \geq 0$, and set
\[ r_{C}(P) = \begin{cases}
            -\infty & \text{if \#$\bigl( B_C (P;0) \bigr) = \infty$,} \\
             \sup \{r \geq 0 \mid \mbox{\#}\bigl( B_C (P;r) \bigr)< \infty \}
              & \text{otherwise.}
\end{cases}   \]
Then, we have the following conjectures due to Bogomolov.

\begin{conj} (Bogomolov's conjecture)\textbf{.}
If $f$ is non-isotrivial, then $r_{C}(P) > 0$ for all $P$.
\end{conj}

\begin{conj} \label{eff-b_conj}
(Effective Bogomolov's conjecture)\textbf{.} If $f$ is non-isotrivial, 
then there exists an effectively calculated positive number $r_0$ such that 
\[ \inf_{P \in \mathrm{Pic}^0(C)(\overline{K})} r_{C}(P)
\geq r_0 . \]
\end{conj}

In order to describe $r_0$ above, 
we introduce the types of nodes of a semistable curve.
Let $C$ be a semistable curve of genus $g$ and $P$ a node of $C$.
We can assign a number $i$ to the node $P$ in the following way.
Let $\nu : C_{P} \to C$ be the partial normalization at $P$.
If $C_{P}$ is connected, then $i = 0$.
Otherwise, $i$ is the minimum of arithmetic genera of two connected components
of $C_P$.
We say the node $P$ of $C$ is of type $i$.
We denote by $\delta_i (C)$
the number of nodes of type $i$,
and by $\delta_i (X/Y)$ the number of nodes of type $i$ in all the fibers
of $f : X \to Y$,
i.e., $\delta_i (X/Y) = \sum_{y \in Y} \delta_i (X_y)$.

Moriwaki proved the following results.
\begin{enumerate}
\renewcommand{\labelenumi}{(\alph{enumi})}
\item (\cite{m-3} etc.) (char$(k)=0$).
If $f$ is not smooth and every singular fiber of $f$ is a tree of 
stable components, then
\[ \inf_{P \in \mathrm{Pic}^0(C)(\overline{K})} r_{C}(P)
\geq 
\sqrt{\frac{(g - 1)^2}{g(2g+1)} \Biggl( \frac{g-1}{3} \delta_0 (X/Y)
+ \sum_{i=1}^{[g/2]} 4i(g-i) \delta_i (X/Y) \Biggl)} \mbox{ . } \]
\item (\cite{m-2}) (char$(k) \geq 0$).
If $g=2$, then $f$ is not smooth and 
\[ \inf_{P \in \mathrm{Pic}^0(C)(\overline{K})} r_{C}(P)
\geq 
\sqrt{\frac{2}{135} \delta_0 (X/Y) + \frac{2}{5} \delta_1 (X/Y)} \mbox{ .} \]
\end{enumerate}

In this paper, we would like to prove the effective Bogomolov's conjecture
for generically smooth semistable hyperelliptic curves.

Let $C$ be a semistable curve over $k$.
We say that $C$ is a semistable hyperelliptic curve if there exist
a valuation ring $R$ with residue field $k$ and a generically smooth
semistable curve $f : Z \to \mathrm{Spec}(R)$ such that
the generic fiber of $f$ is a smooth hyperelliptic curve
and the special fiber of $f$ is $C$.
By the definition, $C$ has an involution $\iota$,
and we can see that $C/ \langle \iota \rangle$ is a nodal curve
which is a tree of $\mathbb{P}^1$.
For details, see \cite{c-h}.

Now let $f : X \to Y$ be a generically 
smooth hyperelliptic semistable curve of genus $g$ with the 
hyperelliptic involution $\iota$.
Let $C$ be a fiber of $f$, which is a semistable hyperelliptic curve over $k$ 
with the involution $\iota = \iota|_{C}$, and
$P$ a node of $C$ of type $0$.
We can also assign a number $j$ to the pair of nodes $(P,\iota(P))$ of type $0$
in the following way.
If $P = \iota(P)$, we set $j=0$. 
If $P \neq \iota(P)$, then  the partial normalization at $P$ and $\iota(P)$ 
$C_{P,\iota(P)}$ has two connected components
since $C/ \langle \iota \rangle$ is a tree of $\mathbb{P}^1$.
We set $j$ to be the minimum of arithmetic genera of two connected components
of $C_{P,\iota(P)}$. 
We say that the node $P$, or the pair of nodes $(P,\iota(P))$ is 
of type $(0,j)$, or of subtype $j$.
We denote by $\xi_{0}(C)$ the number of \emph{nodes} of type $(0,0)$,
and by $\xi_{j}(C)$ the number of \emph{such pairs of nodes} of type $(0,j)$
for $j \geq 1$.
Moreover, we set
\[ \xi_{j} (X/Y) = \sum_{y \in Y} \xi_{j}(X_{y}) . \]

The following are the main results of this paper.

\begin{thm}
\emph{(char$(k) = 0$).}
We assume that $f$ is hyperelliptic. 
Then, Bogomolov's conjecture holds for $f$.
In addition,  $f$ is not smooth and 
\[ \inf_{P \in \mathrm{Pic}^0(C)(\overline{K})} r_{C}(P) \geq \sqrt{r_0} , \]
where $r_0$ is a positive number given below.
\begin{enumerate}
\renewcommand{\labelenumi}{(\arabic{enumi})}
\item
If $g=3,4$, then
\[
\begin{split}
r_0 =&
\frac{(g-1)^2}{g(2g+1)} 
\Biggl(  \frac{(2g-5)}{12} \xi_0(X/Y) \\
&+ \sum_{j=1}^{[(g-1)/2]} (2j(g-1-j)-1) \xi_{j}(X/Y)
+ \sum_{i=1}^{[g/2]} 4i(g-i) \delta_{i}(X/Y)
\Biggr) .
\end{split} \]
\item
If $g \geq 5$, then
\[ \begin{split}
r_0 &=
\frac{(g-1)^2}{g(2g+1)} 
\Biggl(  \frac{(2g-5)}{12} \xi_0(X/Y) \\
&+ \sum_{j=1}^{[(g-1)/2]} \frac{2(3j(g-1-j)-g-2)}{3} \xi_{j}(X/Y)
+ \sum_{i=1}^{[g/2]} 4i(g-i) \delta_{i}(X/Y)
\Biggr).
\end{split}
\]
\end{enumerate}
\end{thm}

Finally, the auther would like to express my deep gratitude to 
Professor Atsushi Moriwaki
for giving me a lot of valuable advice.

\section{Some remarks on the admissible constants}

In this paper, we mean by a graph a topological graph in
sense of \cite{z} equipped with the set of edges and the set of vertices.

Let $G$ be a connected graph,
and $\mathrm{Vert}(G)$ (resp. $\mathrm{Ed}(G)$) the set of vertices
(resp. edges) of $G$.
If $\sim$ is an equivalence relation in $\mathrm{Ed}(G)$, 
we set ${\mathrm{Ed}(G)}^{\sim} = \mathrm{Ed}(G) / \!\!\!\sim$.
Let $\bigoplus_{e \in \mathrm{Ed}(G)} \mathbb{R} e$ 
(resp. $\bigoplus_{\bar{e} \in \mathrm{Ed}(G)^{\sim}} \mathbb{R} \bar{e}$) 
be the $\mathbb{R}$-vector space
formally generated by $\mathrm{Ed}(G)$ (resp. ${\mathrm{Ed}(G)}^{\sim}$),
and $\mathcal{M}(\mathrm{Ed}(G))$ 
(resp. $\mathcal{M}({\mathrm{Ed}(G)}^{\sim})$) the dual vector space of 
$\bigoplus_{e \in \mathrm{Ed}(G)} \mathbb{R} e$
(resp. $\bigoplus_{\bar{e} \in \mathrm{Ed}(G)^{\sim}} \mathbb{R} \bar{e}$).
We express by $\{ {e}^{\ast} \}_{e \in \mathrm{Ed}(G)}$
(resp. $\{ \bar{e}^{\ast} \}_{\bar{e} \in \mathrm{Ed}(G)}$)
the dual basis of $\mathcal{M}(\mathrm{Ed}(G))$
(resp. $\mathcal{M}({\mathrm{Ed}(G)}^{\sim})$)
with respect to
$\mathrm{Ed}(G)$
(resp. ${\mathrm{Ed}(G)}^{\sim}$).
We have the natural projection 
$\bigoplus_{e \in \mathrm{Ed}(G)} \mathbb{R} e 
\twoheadrightarrow 
\bigoplus_{\bar{e} \in \mathrm{Ed}(G)^{\sim}} \mathbb{R} \bar{e}$ 
and the natural inclusion
$\mathcal{M}({\mathrm{Ed}(G)}^{\sim})
 \hookrightarrow \mathcal{M}(\mathrm{Ed}(G))$.
Set
\[ \mathcal{M}({\mathrm{Ed}(G)}^{\sim})_{>0} =
\Biggl\{ \lambda : 
\bigoplus_{\bar{e} \in \mathrm{Ed}(G)^{\sim}} \mathbb{R} \bar{e}
 \to \mathbb{R} \mid
 \lambda(\bar{e}) > 0 \mbox{ for any } e \in \mathrm{Ed}(G) 
\Biggr\}. \]
Note that
to give an element $\lambda \in \mathcal{M}({\mathrm{Ed}(G)}^{\sim})_{>0}$ 
is nothing but
to give length to each edge such that
length of $e$ is $\lambda (\bar{e})$.
In this sense, we sometimes call an element $\lambda \in 
\mathcal{M}({\mathrm{Ed}(G)}^{\sim})_{>0}$ 
a Lebesgue measure on $G$, and call
a graph equipped with a Lebesgue measure a metrized graph.

Now, we recall several facts on Green's function
on a metrized graph. For details on metrized graphs, see \cite{z}.

Let $(G;\lambda)$ be a connected metrized graph and $D$ an 
$\mathbb{R}$-divisor on $G$.
If $\mathrm{deg}(D) \neq -2$, then there are a unique measure 
$\mu_{(G;\lambda ,D)}$ on $G$ and a unique function 
$g_{(G;\lambda,D)}$ on $G \times G$
with the following properties.
\begin{enumerate}
\renewcommand{\labelenumi}{(\alph{enumi})}
\item
$\displaystyle \int_G \mu_{(G;\lambda,D)} = 1$.
\item
$g_{(G;\lambda,D)}(x,y)$ is symmetric and continuous on $G \times G$.
\item
For a fixed $x \in G$, 
$\Delta_{y}(g_{(G;\lambda,D)}(x,y)) = \delta_{x} - \mu_{(G;\lambda,D)}$.
\item
For a fixed $x \in G$, 
$\displaystyle \int_G g_{(G;\lambda,D)}(x,y) \mu_{(G;\lambda,D)}(y) = 0$.
\item
$g_{(G;\lambda,D)}(D,y) + g_{(G;\lambda,D)}(y,y)$ 
is a constant function on $y \in G$.
\end{enumerate}
The constant $g_{(G;\lambda,D)}(D,y) + g_{(G;\lambda,D)}(y,y)$ 
is denoted by $c(G;\lambda,D)$.
Further, we set
\[ \epsilon(G;\lambda,D) =
 2 \mathrm{deg}(D) c(G;\lambda,D) - g_{(G;\lambda,D)}(D,D), \]
which we call the \emph{admissible constant} of $(G;\lambda,D)$.
In this paper, we consider polarizations on $G$ supported in 
$\mathrm{Vert}(G)$ only.

Let $(G;\lambda,D)$ be a connected polarized metrized graph.
We consider the following constants arising from $(G;\lambda,D)$.
(In the following, 
$e$ is an edge, $P_e$ and $Q_e$ are the terminal points of $e$,
$e^{\circ} = e \setminus \{ P_e,Q_e \}$, and $P,Q \in \mathrm{Vert}(G)$.)
\begin{align*}
&l_{\bar{e}}= \lambda(\bar{e}) &\quad&\text{: the length of $e$}\\
&r_{(G;\lambda)} (P,Q) &\quad&\text{: the resistance between $P$ and $Q$}\\
&r_{e} = r_{G \setminus  \{e^{\circ}\}} (P_e,Q_e) \\
&g_{(G;\lambda,D)}(P,Q)  \\
&\epsilon(G;\lambda,D)
\end{align*}
If $\gamma(G;\lambda,D)$ is one of the above constants,
then it is easy to see that the function on 
$\mathcal{M}({\mathrm{Ed}(G)}^{\sim})_{>0}$ defined by
\[ \lambda \mapsto \gamma (G;\lambda,D) \]
is a rational function.
We denote these functions by the ``similar'' symbols as follows.
\begin{align*}
&X_{\bar{e}} : \lambda \mapsto \lambda(\bar{e}) \\
&r_{G} (P,Q) : \lambda \mapsto r_{(G;\lambda)} (P,Q) \\
&R_{e} : \lambda \mapsto r_{e} \\
&g_{(G,D)}(P,Q) : \lambda \mapsto  g_{(G;\lambda,D)}(P,Q)\\
&\epsilon(G,D) : \lambda \mapsto \epsilon(G;\lambda,D)
\end{align*}
When we do not have to emphasize $\lambda$,
we sometimes write $\bar{\epsilon}(G,D)$ for $\epsilon(G;\lambda,D)$,
for example.
Note that these rational functions can be viewed as elements of rational
function field
$\mathbb{Q}\bigl(\{ X_{\bar{e}} \}_{\bar{e} \in {\mathrm{Ed}(G)}^{\sim}}\bigr)$
generated by indeterminates 
$\{ X_{\bar{e}} \}_{\bar{e} \in {\mathrm{Ed}(G)}^{\sim}}$.

Let $G$ be a connected graph and $S$ a subset of 
${\mathrm{Ed}(G)}^{\sim}$.
We define $G_{S}$ as the graph obtained by contracting all edges $e$
with $\bar{e} \in S$, and define
$G^S$ as $G_{\mathrm{Ed}(G)^{\sim} \setminus S}$.
For a polarization $D$ on $G$, we also define $D_{S}$ (resp. $D^{S}$)
as the polarization on $G_{S}$ (resp. $G^S$) induced by $D$
in the following way.
Let $v$ be a vertex of $G_S$ and $\{ v_{1}, \ldots , v_{k} \}$ the set of
vertices of G which go to $v$ when we contract the edges in $S$.
Then, we set the coefficient of $v$ of $D_S$ to be the sum of coefficients of
all $v_{i}$'s of $D$.
Note that $\deg (D_S) = \deg(D)$.

\begin{lem}
In the same notation as above, we have
\[ \epsilon(G,D)(X_{\bar{e}_0}=0) = \epsilon(G_{\bar{e}_0}, D_{\bar{e}_0}) \]
for any $\bar{e}_0 \in {\mathrm{Ed}(G)}^{\sim}$.
\end{lem}

\begin{Pf}
It is sufficient to show that
\[ \lim_{l_{\bar{e}_0} \to 0} \epsilon(G;\lambda,D) =
\epsilon(G_{\bar{e}_0};\lambda', D_{\bar{e}_0}) \]
for any $\lambda = \sum_{\bar{e} \in {\mathrm{Ed}(G)}^{\sim}}
l_{\bar{e}} \bar{e}^{\ast}$,
where $\lambda' = \sum_{\bar{e} \in {\mathrm{Ed}(G_{\bar{e}_0})}^{\sim}}
l_{\bar{e}} \bar{e}^{\ast}$.

We may assume that all edges are connected closed interval.
Let $s_{e} : e \to [0,l_e]$ be a parameterization.
Then, we can set
\begin{align*} &g(x) := g_{(G;\lambda,D)}(O,x)=
 \alpha_{e} s_{e}(x)^2 + \beta_{e} s_{e}(x) + \gamma_e 
&\quad & \text{on $e$ of $G$} \\
&g'(x) := g_{(G_{\bar{e}_{0}};\lambda',D_{0})}(O,x)=
 \alpha_{e}' s_{e}(x)^2 + \beta_{e}' s_{e}(x) + \gamma_e' 
&\quad & \text{on $e$ of $G_{\bar{e}_{0}}$}
\end{align*}
for some $\alpha_{e},\beta_{e},\gamma_e \in \mathbb{R}$
and $\alpha_{e}',\beta_{e}',\gamma_e' \in \mathbb{R}$.

The continuous condition on $G$, $\Delta(g)=\delta_O - \mu_{(G;\lambda,D)}$,
the continuous condition on $G'$ and 
$\Delta(g')=\delta_O - \mu_{(G_{\bar{e}};\lambda',D_{\bar{e}})}$
give a system of linear equations on $\alpha_{e},\beta_{e}$
and $\alpha_{e}',\beta_{e}'$.
It is easy to see that
$\alpha_{e} \to \alpha_{e}'$, $\beta_{e} \to \beta_{e}'$ when
$l_{\bar{e}} \to 0$.
By the conditions $\int_{G} g \mu_{(G;\lambda,D)}$
and $\int_{G_{\bar{e}}} g \mu_{(G_{\bar{e}};\lambda',D_{\bar{e_{0}}})}$,
we also have $\gamma_{e} \to \gamma_{e}'$,
hence we obtain the lemma.
\end{Pf}

Let $G_{1}$ and $G_{2}$ be graphs.
Fix vertices $v_{1} \in G_{1}$ and $v_{2} \in G_{2}$.
The \emph{one-point-sum} $G_{1} \vee G_{2}$ with respect to $v_{1}$
and $v_{2}$ is defined as $(G_{1} \amalg G_{1})/v_{1}\sim v_{2}$.
The set of edges is naturally defined by 
$\mathrm{Ed}(G_{1} \vee G_{2}) = \mathrm{Ed}(G_{1}) \amalg \mathrm{Ed}(G_{2})$
and the set of vertices is defined by
$\mathrm{Vert}(G_{1} \vee G_{2})
= (\mathrm{Vert}(G_{1}) \amalg \mathrm{Vert}(G_{2}))/v_{1}\sim v_{2}$.
If $G_{i}$ has a Lebesgue measure $\lambda_{i}$ for $i = 1,2$,
then $G_{1} \vee G_{2}$ has the canonical Lebesgue measure given by
$\lambda = \lambda_{1} + \lambda_{2}$.

\begin{Def}
Let $G$ be a connected graph which is not one point.
$G$ is said to be \emph{reducible}
if there exist two graphs $G_{1}$ and $G_{2}$ which are not one point
such that $G$ is a one-point-sum of $G_{1}$ and $G_{2}$.
$G$ is said to be \emph{irreducible} if it is not reducible.
\end{Def}

For any connected graph $G$,
we have the \emph{irreducible decomposition} of $G$.
Set
\[ J = \big\{ v \in \mathrm{Vert}(G) \mid
G \setminus \{ v \} 
\mbox{ is not connected.} \big\} . \]
Let $H^{\circ}$ be a connected component of $G \setminus J$.
Then, the closure of $H^{\circ}$ is an irreducible subgraph of $G$, and
we call it an irreducible component of $G$.
Let $\{ G_1, \ldots, G_n \}$ be the set of irreducible components.
For a permutation $\sigma : (1,2, \ldots, n) \mapsto 
(i_1, i_2, \ldots, i_n)$, we define a sequence of subgraphs of $G$ 
inductively by
\[ H_1^{\sigma} = G_{i_1}, \quad 
H_{k}^{\sigma} = H_{k-1}^{\sigma} \cup G_{i_{k}} .\]
Then, we can easily see by the definition of $\{ G_1, \ldots, G_n \}$ that
there exists a permutation $\sigma$ such that
$H_{k-1}^{\sigma} \cup G_{i_{k}}$ is a one-point-sum 
of $H_{k-1}^{\sigma}$ and $G_{i_{k}}$ for all $k = 2, \ldots,n$.
In this sense, we write 
$G_1 \vee\cdots\vee G_n$ instead of $G_1 \cup \cdots \cup G_n$
and call it the \emph{irreducible decomposition}
of $G$.
We denote by $\mathrm{Irr}(G)$ the set of all irreducible components of $G$.

The next proposition implies that irreducible graphs are fundamental
for calculating the admissible constants.

\begin{prop} \label{sum-formula}
Let $G_1$, $G_2$ and $G=G_1 \vee G_2$ be connected graphs,
and $D$ a polarization supported in $\mathrm{Vert}(G)$
with $\deg D \neq -2$.
Let $D_i$ be the polarization on $G_i$ defined by 
$D_i = D^{\mathrm{Ed}(G_i)}$ for $i=1,2$.
Then, we have
\[ \epsilon(G,D) = \epsilon(G_1,D_1) + \epsilon(G_2,D_2) \]
as rational functions on 
$\mathcal{M}({\mathrm{Ed}(G)})_{>0}$.
\end{prop}

\begin{Pf}
Let $\lambda_1$ and $\lambda_2$ be Lebesgue measures on $G_1$ and $G_2$
respectively.
$\lambda = \lambda_1 + \lambda_2$ is a Lebesgue measure on $G$.
By \cite [Lemma 3.7] {z}, we have
\[ \mu_{(G;\lambda,D)} 
= \mu_{(G_1;\lambda_{1},D_1)} + \mu_{(G_2;\lambda_{2},D_2)} - \delta_O ,\]
where $\{ O \} = G_{1} \cap G_{2}$.
Consider the following function on $G$:
\[
g(x)=
\begin{cases}
g_{(G_1;\lambda_{1},D_1)} ( O, x) + g_{(G_2;\lambda_{2},D_2)} ( O, O)  
& \text{if $x \in G_1$,} \\
g_{(G_2;\lambda_{2},D_2)} ( O, x) + g_{(G_1;\lambda_{1},D_1)} ( O, O)  
& \text{if $x \in G_2$.}
\end{cases}
\]
Then, we can easily check that $g$ is continuous on $G$,
$\Delta(g) = \delta_O - \mu_{(G;\lambda,D)}$, and 
$\int_{G} g \mu_{(G;\lambda,D)} = 0$.
Thus we have $g_{(G;\lambda,D)}(O,x) = g(x)$.
Therefore, by \cite [Lemma 4.1] {m-4},
we obtain the formula.
\end{Pf}

\section{Calculations of the admissible constants for hyperelliptic graphs}

\subsection{Definitions and terminology}

First of all, we give the definition of a particular class of graphs,
called hyperelliptic graphs.

\begin{Def} \label{def. of H.G}
Let $G$ be a connected graph, and
$\mathrm{Vert}(G)$ (resp. $\mathrm{Ed}(G)$) the set of vertices
(resp. edges) of G.
Suppose that $G$ has a homeomorphism $\iota : G \to G$
such that $\iota^2$ is the identity on $G$, called the 
involution on $G$, which induces naturally an automorphism on 
$\mathrm{Vert}(G)$
and $\mathrm{Ed}(G)$ respectively.
Then, $(G,\mathrm{Vert}(G),\mathrm{Ed}(G),\iota)$, or simply $G$, is called
a \emph{hyperelliptic graph} if it has the following properties.
\begin{enumerate}
\renewcommand{\labelenumi}{(\arabic{enumi})}
\item Every edge is homeomorphic to the connected closed interval.
\item $\iota(e) \neq e$ for any $e \in \mathrm{Ed}(G)$.
\item If $v$ is a vertex with $\iota(v) \neq v$, then there exist
   at least three edges which start from $v$.
\item The topological space $G/ \langle \iota \rangle$ has no
loops. (We call such a graph a \emph{tree}.)
\end{enumerate}
\end{Def}

Note that $G/ \langle \iota \rangle$ is a connected graph 
whose vertices and edges are given by 
${\mathrm{Vert}(G)}^{\sim} = \mathrm{Vert}(G) / \langle \iota \rangle$
and ${\mathrm{Ed}(G)}^{\sim} = \mathrm{Ed}(G) / \langle \iota \rangle$
respectively.

When we talk on a measure on a hyperelliptic graph,
we always assume that it is invariant under the involution, i.e.,
an element of $\mathcal{M}({\mathrm{Ed}(G)}^{\sim})_{>0}$ 
with respect to the equivalence
relation arising from $\iota$.

\begin{ex} \label{ex.of.h.g}
We shall give an example of hyperelliptic graphs, which is the main object
in this paper.
Let $G_{1}$ be the metrized graph by the configuration of
a singular fiber $C_{1}$ of semistable hyperelliptic curve
$f : X \to Y$ as in the introduction.
We assume that $C_{1}$ does not have nodes of positive type.
$\mathrm{Vert}(G_{1})$ and $\mathrm{Ed}(G_{1})$ correspond to 
the set of irreducible components of $C_{1}$ and the set of 
nodes of $C_{1}$ respectively.
The hyperellptic involution $\iota$ also acts on $\mathrm{Vert}(G_{1})$
and $\mathrm{Ed}(G_{1})$.
Then, there may exists an edge $e$ with $\iota(e) = e$.
Note that if such $e$ is the connected closed inteval,
then the vertices which are the terminal points of $e$
are moved to each other by $\iota$.
For $e$ with $\iota(e) = e$, let $v_{e}$ be the point on $e$ such that
$e \setminus ( \{ \mbox{vertices on $e$} \} \cup \{ v_{e} \} )$
is a disjoint union of two open segments of same length.
Now, let $G_{2}$ be the metrized graph
which is same as $G_{1}$ as a metrized topological space, such that
$\mathrm{Vert}(G_{2})$ is the union of $\mathrm{Vert}(G_{1})$ and 
the set of such $v_{e}$'s as above,
and that $\mathrm{Ed}(G_{2})$ is the segments in $G_{1}$ which connect
two points in $\mathrm{Vert}(G_{2})$.
Then, we can make $\iota$ act on $G_{2}$ such that
$\iota$ is a symmetric homeomorphism
and $\iota(e) \neq e$ for any $e \in \mathrm{Ed}(G_{2})$.
Let $V$ be the subset of $\mathrm{Vert}(G_{2})$ consisting of 
vertices $v$ such that $\iota(v) \neq v$ and there are
only two edges starting from $v$.
Then, $G_{2} \setminus \bigl(\mathrm{Vert}(G_{2}) \setminus V\bigr)$
is a disjoint union of open segments.
Let $G_{3}$ be the metrized graph which is nothing but $G_{2} = G_{1}$
as a metrized space,
such that $\mathrm{Vert}(G_{3}) = \mathrm{Vert}(G_{2}) \setminus V$
and $\mathrm{Ed}(G_{3})$ is the set of segments in $G_{2}$ which connect
two points in $\mathrm{Vert}(G_{3})$.
We can also make $\iota$ act on $G_{3}$ naturally such that
$\iota$ is a symmetric homeomorphism.
Noting, in addition,  that $C_{1} / \langle \iota \rangle$ is a tree of 
$\mathbb{P}^{1}$,
we can easily see by its construction that
$G_{3}$ is a hyperelliptic graph with an $\iota$-invariant measure.
\end{ex}

We fix the following terminology.

\begin{Def}
Let $G$ be a hyperelliptic graph.
\begin{enumerate}
\renewcommand{\labelenumi}{(\arabic{enumi})}
\item $v \in \mathrm{Vert}(G)$ is said to be \emph{fixed} if $\iota(v)=v$.
   We denote by $\mathrm{Vert}_{\mathrm{f}}(G)$ the set of fixed vertices.
\item $v \in \mathrm{Vert}(G)$ is said to be \emph{non-fixed}
   if $\iota(v) \neq v$.
   We denote by $\mathrm{Vert}_{\mathrm{n.f}}(G)$ 
   the set of non-fixed vertices.
\item $e \in \mathrm{Ed}(G)$ is said to be \emph{disjoint} if 
   $e \cap \iota(e) = \emptyset$.
   We denote by $\mathrm{Ed}_{0}(G)$ the set of disjoint edges.
\item $e \in \mathrm{Ed}(G)$ is said to be \emph{one-jointed} if 
   $e \cap \iota(e)$ is a set of one point. 
   We denote by $\mathrm{Ed}_{1}(G)$ the set of one-jointed edges.
\item $e \in \mathrm{Ed}(G)$ is said to be \emph{two-jointed} if 
   $e \cap \iota(e)$ is a set of two points.
   We denote by $\mathrm{Ed}_{2}(G)$ the set of two-jointed edges.
\end{enumerate}
\end{Def}
If $S$ is one of the above sets, we denote by $S^{\sim}$
the set $S/\langle \iota \rangle$,
and we write $\bar{s}$ for the class of $s \in S$ in $S^{\sim}$.

Let us consider several lemmas concerning the above definitions.
\begin{lem}
If $G_1$ is an irreducible component of a hyperelliptic graph $G$,
then we have $\iota(G_1) = G_1$.
\end{lem}

\begin{Pf}
Let $\pi : G \to G / \langle \iota \rangle$ be the natural projection.
Suppose $\iota(G_1) \neq G_1$.
Let $P$ be a vertex of $G_1$ which joints $G_1$ with another component,
and $H$ the subgraph containing $G_1$ such that 
$H \setminus \{ P , \iota(P) \}$ is the connected component of 
$G \setminus \{ P , \iota(P) \}$.
We note the following claim.
\begin{cl}
Let $G_{2}$ be an irreducible component of $G$ 
different from $G_{1}$ having $P$ as a jointing point
with $G_{1}$. Then, $G_{2} \cap H = \{ P \}$.
\end{cl}
\begin{Pf}
Since $P$ is a jointing point,
$G_{2} \setminus \{ P \}$ is not contained in the connected component
of $G \setminus \{ P \}$ which $G_{1} \setminus \{ P \}$ belongs to.
\end{Pf}

Now, assume $P = \iota(P)$.
Then, $\iota(G_{1})$ is an irreducible component of $G$ with 
$P \in \iota(G_{1})$.
Since $\iota(G_{1}) \neq G_{1}$,
we must have $\iota(G_{1}) \cap H = \{ P \}$ by the claim,
hence $\iota(H) \neq H$.
Therefore, we see $\iota(H) \cap H = \{ P \}$ by the definition of $H$,
and $\pi |_{H} : H \to \pi(H)$ is an isomorphism, accordingly,
$H$ is a tree.
Take a terminal point $Q$ of $H$ different from $P$.
Then, since $H$ is a tree,
there exists only one edge starting from $Q$.
This contradicts to the definition of hyperelliptic graphs.
Therefore, we must have $\iota(P) \neq P$.
Noting $\iota(H) \cap H \supset \{ P , \iota(P) \}$ and the definition
of $H$, we see again $H = \iota(H)$.
Since $P$ and $\iota(P)$ are two distinct jointing points,
$G \setminus (H \setminus \{ P , \iota(P) \} )$
has two connected components $G_{3}$ and $G_{4}$
with $\iota(G_{3}) = G_{4}$.
Therefore,
$\pi |_{G_{3}} : G_{3} \to \pi(G_{3})$ is an isomorphism,
and we have a contradiction in the same way.
\end{Pf}

\begin{lem}
If $G = G_1 \vee G_2$ is a hyperelliptic graph,
then $\iota(O)=O$, where $O$ is the jointing point of 
$G_1$ and $G_2$ in $G$.
\end{lem}

\begin{Pf}
If $\iota(O) \neq O$, then $G_1 \cap G_2$ has two points 
$O$ and $\iota(O)$ by the above lemma, which contradicts to the assumption
of this lemma.
\end{Pf}
In virtue of the above two lemmas, we see that $G$ is hyperelliptic
if and only if every irreducible component is hyperelliptic.

The next lemma characterizes jointing points of a hyperelliptic graph.
\begin{lem} \label{jointing-points-lemma}
Let $G$ be a hyperelliptic graph, and $v$ a vertex of $G$.
Then, $v$ is a jointing point of irreducible components if and only if
$v$ is fixed and at least four edges start from $v$.
\end{lem}

\begin{Pf}
The ``only if'' part is obvious from the above rwo lemmas.
We will show the ``if'' part.

Let $e_{1}$,  $\iota(e_{1})$, $e_{2}$ and $\iota(e_{2})$
be four edges starting from $v$.
If $v$ is not a jointing point, then $G \setminus \{ v \}$ is connected,
hence we can find a path connecting $e_{1}$ with $e_{2}$ in 
$G \setminus \{ v \}$.
This shows, however, that $G / \langle \iota \rangle$ has a loop,
which contradicts to Definition
\ref{def. of H.G} (4).
\end{Pf}

\begin{lem}
If $e$ is a two-jointed edge, 
then $e \cup \iota(e)$ is an irreducible component
of $G$.
\end{lem}

\begin{Pf}
By Lemma \ref{jointing-points-lemma},
it is sufficient to show that
if $v$ is a vertex of $G$ with $\iota(v) \neq v$,
then $\iota(v)$ and $v$ cannot be connected by one edge.
Suppose that $\iota(v)$ and $v$ are connected by one edge $e$.
If we suitably parameterize $e \cup \iota(e)$,
then it is homeomorphic to the circle $S^1 =\{(\cos t ,\sin t) \}$ 
and the action of $\iota$ on $e \cup \iota(e)$ is nothing but a
map from $S^1$ to $S^1$ given by $t \mapsto t+\pi$.
Therefore, the image of $e \cup \iota(e)$ in $G / \langle \iota \rangle$ 
is a circle,
which is a contradiction.
\end{Pf}
Thus, we can give the following definitions.

\begin{Def}
\begin{enumerate}
\renewcommand{\labelenumi}{(\arabic{enumi})}
\item An irreducible hyperelliptic graph $G$ is said to be \emph{simple} if
$G$ consists of two two-jointed edges.
\item A hyperelliptic graph $G$ is said to be \emph{semisimple} if
every irreducible component of $G$ is simple.
\end{enumerate}
\end{Def}
We can easily see that the simple graph is uniquely determined.
We denote by $SG$ the simple graph.
(See Figure \ref{simplegraph}.)

\begin{figure}[htbp]
\begin{center}
\begin{picture}(100,80)
\put(10,40){\circle*{5}}
\put(90,40){\circle*{5}}
\qbezier(10,40)(50,70)(90,40)
\qbezier(10,40)(50,10)(90,40)
\end{picture}
\caption{$SG$} \label{simplegraph}
\end{center}
\end{figure}

Let $\mbox{\#} \mathrm{Irr}(G)$ be the number of irreducible components of $G$
and let $\mbox{\#} \mathrm{Irr}_\mathrm{s}(G)$ be
that of irreducible components of $G$
which is simple.

\begin{rem}
The following are immediate from the definitions.
\begin{enumerate}
\renewcommand{\labelenumi}{(\arabic{enumi})}
\item We have $\mbox{\#} \mathrm{Irr}(G)=\mbox{\#} \mathrm{Irr}(G_{\bar{e}})$
for $\bar{e} \in {\mathrm{Ed}_{0}(G)}^{\sim}$.
\item We have $\mbox{\#} \mathrm{Irr}(G)<\mbox{\#} \mathrm{Irr}(G_{\bar{e}})$
for $\bar{e} \in {\mathrm{Ed}_{1}(G)}^{\sim}$.
\item 
We have $\mbox{\#} \mathrm{Irr}(G)=\mbox{\#} \mathrm{Irr}(G_{\bar{e}}) +1$
for $\bar{e} \in {\mathrm{Ed}_{2}(G)}^{\sim}$.
\end{enumerate}
\end{rem}

\begin{Def}
Let $G$ be a hyperelliptic graph.
We define the \emph{size} of $G$, denoted by $\mathrm{sz}(G)$,
in the following way.
\begin{enumerate}
\renewcommand{\labelenumi}{(\arabic{enumi})}
\item If $G$ is irreducible, we set
\[ \mathrm{sz}(G) = 
     \begin{cases}
1 & \mbox{if $G$ is simple,} \\
\mathrm{\#} ({\mathrm{Ed}_{1}(G)}^{\sim}) - 1 & \mbox{otherwise.}
\end{cases} \]
\item 
If $G = G_1 \vee G_2 \vee \cdots \vee G_k$, where $G_i$ ($i= 1,2,\cdots,k$)
is irreducible, we set
\[ \mathrm{sz}(G) = \sum_{i=1}^{k} \mathrm{sz}(G_i). \]
\end{enumerate}
\end{Def}
Note that
\begin{align*} 
\mathrm{sz}(G)
&= 
\mathrm{\#} ({\mathrm{Ed}_{1}(G)}^{\sim}) 
- ( \mbox{\#} \mathrm{Irr}(G) - \mbox{\#} \mathrm{Irr}_{\mathrm{s}}(G) )
+ \mbox{\#} \mathrm{Irr}_{\mathrm{s}}(G) \\
&= \mathrm{\#} ({\mathrm{Ed}_{1}(G)}^{\sim}) - \mbox{\#} \mathrm{Irr}(G)
+ 2 \mbox{\#} \mathrm{Irr}_{\mathrm{s}}(G).
\end{align*}

\begin{figure}[htbp]
\begin{center}
\begin{picture}(120,120)
\put(60,110){\circle*{5}}
\put(60,10){\circle*{5}}
\put(10,60){\circle*{5}}
\put(110,60){\circle*{5}}
\put(20,60){\circle*{5}}
\put(30,60){\circle*{5}}
\put(100,60){\circle*{5}}
\put(10,60){\line(1,1){50}}
\put(10,60){\line(1,-1){50}}
\put(20,60){\line(4,5){40}}
\put(20,60){\line(4,-5){40}}
\put(30,60){\line(3,5){30}}
\put(30,60){\line(3,-5){30}}
\put(100,60){\line(-4,5){40}}
\put(100,60){\line(-4,-5){40}}
\put(110,60){\line(-1,1){50}}
\put(110,60){\line(-1,-1){50}}
\put(50,57){$\cdots \cdots$}
\put(60,110){\line(-2,-5){10}}
\put(60,10){\line(-2,5){10}}
\put(60,110){\line(3,-5){15}}
\put(60,10){\line(3,5){15}}
\end{picture}
\caption{$\mathbf{G}_n$} \label{elementarygraph}
\end{center}
\end{figure}
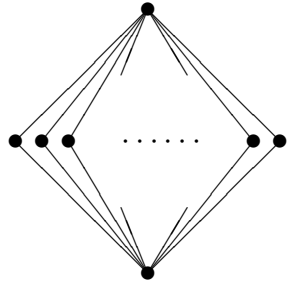

\begin{Def}
An irreducible hyperelliptic graph $G$ of size $n$ is said to be 
(\emph{$n$-th}) \emph{elementary} if all edges of $G$ are one-jointed.
\end{Def}
Note that $\mathrm{sz}(G) >1$ if $G$ is an elementary graph.
The $n$-th elementary graph is uniquely determined, and
we denote it by $\mathbf{G}_n$. (See Figure \ref{elementarygraph}.)

\begin{lem}
Let $G$ be a hyperelliptic graph.
Then, we have $\mathrm{sz}(G) = \mathrm{sz}(G_{\bar{e}})$
for $e \in \mathrm{Ed}_{0}(G) \cup \mathrm{Ed}_{1}(G)$. 
\end{lem}

\begin{Pf}

We may assume that $G$ is irreducible.

If $e$ is disjoint, then it is obvious.
Suppose $e$ is one-jointed.
Let $v$ and $\iota(v)$ be the terminal points of $e \cup \iota(e)$,
and suppose that $k$ disjoint edges and $(l+1)$ one-jointed edges
start from $v$.
Then, $G_{\bar{e}}$ decomposes into $k$ non-simple components
and $l$ simple components.
Hence, we have
\begin{align*}
\mathrm{sz}(G_{\bar{e}}) &=
{\mathrm{\#}\mathrm{Ed}_{1}(G_{\bar{e}})}^{\sim} 
-k+l \\
&= \bigl( \mathrm{\#}{\mathrm{Ed}_{1}(G)}^{\sim} -1-l+k \bigr) -k+l \\
&= \mathrm{sz}(G).
\end{align*}
\end{Pf}

We define functions $\nu_{G}^{0}$, $\nu_{G}^{1}$, and $\nu_G$ 
on $\mathrm{Vert}_{\mathrm{n.f}}(G)$ as follows.
\[ \begin{split}
&\nu_{G}^{0} (v) = \text{the number of disjoint edges which start from $v$.} \\
&\nu_{G}^{1} (v) 
= \text{the number of one-jointed edges which start from $v$.} \\
&\nu_G (v) = \nu_{G}^{0} (v)+\nu_{G}^{1} (v).
\end{split} \]

Let $G$ be a hyperelliptic graph.
Let $\{ X_{\bar{e}} \}_{\bar{e} \in {\mathrm{Ed}(G)}^{\sim}}$ 
be the set of symbols as in section 2,
and $V$ the $\mathbb{Q}$-vector space with basis 
$\{ X_{\bar{e}} \}_{\bar{e} \in {\mathrm{Ed}(G)}^{\sim}}$.
We denote the space of homogeneous polynomials of degree $d$
by $S^d V$, i.e.,
$S^d V$ is the $d$-th symmetric tensor product of $V$.

Put $n = \mathrm{sz}(G)$.
Now, we introduce important polynomials.
\begin{enumerate}
\renewcommand{\labelenumi}{(\arabic{enumi})}
\item Choose distinct 
$\bar{e}_1, \ldots , \bar{e}_n \in {\mathrm{Ed}(G)}^{\sim}$.
Set 
\[ 
\delta_{\bar{e}_1, \ldots , \bar{e}_n} = 
\left\{ \begin{array}{@{\,}ll}
1 & \mbox{if $G^{\bar{e}_1, \ldots , \bar{e}_n}$ is a semisimple hyperelliptic
graph of size $n$,} \\
0 & \mbox{otherwise,}
\end{array} \right. \]
and
\[ L_G = \sum_{
\begin{smallmatrix}
\bar{e}_1, \ldots , \bar{e}_n \\
\mbox{\scriptsize{all distinct}}
\end{smallmatrix}} 
\delta_{\bar{e}_1, \ldots , \bar{e}_n} X_{\bar{e}_1} \cdots X_{\bar{e}_n}
\quad \in S^n V. \]
\item 
Choose distinct 
$\bar{e}_1, \ldots , \bar{e}_{n+1} \in {\mathrm{Ed}(G)}^{\sim}$.
If \#$({\mathrm{Vert}_{\mathrm{n.f}}(G^{\bar{e}_1, \ldots , \bar{e}_{n+1}})}^{\sim})=1$,
we denote by $\nu(\bar{v})$ the number of edges which start from 
a representative $v$ of the unique
non-fixed vertex class $\bar{v}$ of $G^{\bar{e}_1, \ldots , \bar{e}_{n+1}}$.
Then we set
\[ c_{\bar{e}_1, \ldots , \bar{e}_{n+1}} =
\begin{cases}
\nu(\bar{v}) - 2  & \mbox{if
\#$({\mathrm{Vert}_{\mathrm{n.f}}(G^{\bar{e}_1, \ldots , \bar{e}_{n+1}})}^{\sim})=1$} \\
0 & \mbox{otherwise,}
\end{cases} \]
and 
\[ M_G = \sum_{
\begin{smallmatrix}
\bar{e}_1, \ldots , \bar{e}_{n+1} \\
\mbox{\scriptsize{all distinct}}
\end{smallmatrix}}
c_{\bar{e}_1, \ldots , \bar{e}_{n+1}} X_{\bar{e}_1} \cdots X_{\bar{e}_{n+1}}
\quad \in S^{n+1} V. \]
Note that $M_G = 0$ if $G$ is semisimple.
\end{enumerate}

In the case where $G$ is irreducible, 
$L_G$ and $M_G$ can be expressed in another way.
Let $\bar{e}_1, \ldots, \bar{e}_k$ be distinct disjoint edges,
and $v$ a non-fixed vertex of the graph 
$G_{{\mathrm{Ed}_{0}(G)}^{\sim} 
\setminus \{ \bar{e}_1, \ldots, \bar{e}_k \}}$.
For simplicity, we set 
$G':=G_{{\mathrm{Ed}_{0}(G)}^{\sim} 
\setminus \{ \bar{e}_1, \ldots, \bar{e}_k \}}$.
Set $\sigma_{(\bar{e}_1, \ldots, \bar{e}_k; \bar{v})}$ to be the
$(\nu_{G'}^{1}(v)-1)$-th elementary symmetric polynomial on
$\{ X_{\bar{e}} \}_{\bar{e} \in {{\mathrm{Ed}_{1}(G')}^{\sim}_{\bar{v}}}}$
and $\tau_{(\bar{e}_1, \ldots, \bar{e}_k; \bar{v})}$ to be the
$\nu_{G'}^{1}(v)$-th elementary symmetric polynomial on
$\{ X_{\bar{e}} \}_{\bar{e} \in {{\mathrm{Ed}_{1}(G')}^{\sim}_{\bar{v}}}}$,
where ${\mathrm{Ed}_{1}(G')}^{\sim}_{\bar{v}}$ is the set of edge classes
whose representatives start from $v$.

Then, we can easily see that
\begin{align*}
&L_{G}=
\sum_{
\begin{smallmatrix}
\bar{e}_1, \ldots, \bar{e}_k \in  {\mathrm{Ed}_{0}(G)}^{\sim} \\
\mbox{\scriptsize {all distinct}}
\end{smallmatrix}}
\Biggl( \prod_{\bar{v}}
\sigma_{(\bar{e}_1, \ldots, \bar{e}_k; \bar{v})} \Biggr)
X_{\bar{e}_1} \cdots X_{\bar{e}_k} \\
&M_{G}=
\sum_{
\begin{smallmatrix}
\bar{e}_1, \ldots, \bar{e}_k \in  {\mathrm{Ed}_{0}(G)}^{\sim} \\
\mbox{\scriptsize {all distinct}}
\end{smallmatrix}} 
\Biggl(
\sum_{\bar{v}} \Biggl( ( \nu(\bar{v}) - 2 )
\tau_{(\bar{e}_1, \ldots, \bar{e}_k; \bar{v})}
\prod_{\bar{v}' \neq \bar{v}}
\sigma_{(\bar{e}_1, \ldots, \bar{e}_k; \bar{v}')} \Biggr) \Biggr)
X_{\bar{e}_1} \cdots X_{\bar{e}_k},
\end{align*}
where $\bar{v}$ runs over 
${\mathrm{Vert}_{\mathrm{n.f}}(G_{{\mathrm{Ed}_{0}(G)}^{\sim} 
\setminus \{ \bar{e}_1, \ldots, \bar{e}_k \}})}^{\sim}$.
Note that $L_G$ is an irreducible polynomial if $G$ is an
irreducible hyperelliptic graph.

\begin{rem}
Let $G = G_1 \vee G_2$ be a hyperelliptic graph of size $n$.
\begin{enumerate}
\renewcommand{\labelenumi}{(\arabic{enumi})}
\item 
For distinct
$\bar{e}_1, \ldots , \bar{e}_k \in {\mathrm{Ed}(G_1)}^{\sim}$ and distinct
$\bar{e}_{k+1}, \ldots , \bar{e}_n \in {\mathrm{Ed}(G_2)}^{\sim}$,
$G^{\bar{e}_1, \ldots , \bar{e}_n}$ is a semisimple hyperelliptic graph
of size $n$ if and only if 
$G_{1}^{\bar{e}_1, \ldots , \bar{e}_k}$ 
(resp. $G_{2}^{\bar{e}_{k+1}, \ldots , \bar{e}_n}$)
is a semisimple hyperelliptic graph
of size $\mathrm{sz}(G_1)$ (resp. $\mathrm{sz}(G_2)$).
\item
For distinct
$\bar{e}_1, \ldots , \bar{e}_k \in {\mathrm{Ed}(G_1)}^{\sim}$ and distinct
$\bar{e}_{k+1}, \ldots , \bar{e}_{n+1} \in {\mathrm{Ed}(G_2)}^{\sim}$,
$G^{\bar{e}_1, \ldots , \bar{e}_{n+1}}$ is a one-point-sum of 
semisimple hyperelliptic graphs and the $l$-th elementary hyperelliptic graph
if and only if one of $\{ G_{1}^{\bar{e}_1, \ldots , \bar{e}_k},
G_{2}^{\bar{e}_{k+1}, \ldots , \bar{e}_{n+1}} \}$, say 
$G_{1}^{\bar{e}_1, \ldots , \bar{e}_k}$,
is a semisimple hyperelliptic graph of size $\mathrm{sz}(G_1)$,
and the other, say $G_{2}^{\bar{e}_{k+1}, \ldots , \bar{e}_{n+1}}$, is 
a one-point-sum of $\bigl(\mathrm{sz}(G_{2})-l \bigr)$ 
simple hyperelliptic graphs 
and the $l$-th elementary hyperelliptic graph.
\end{enumerate}
\end{rem}
The next lemma is simple, but important for our latter purpose.
\begin{lem} \label{L-M-lem}
Let $G$ be a hyperelliptic graph.
\begin{enumerate}
\renewcommand{\labelenumi}{(\arabic{enumi})}
\item
$L_{G_{\bar{e}}} = L_{G}(X_{\bar{e}}=0)$,
$M_{G_{\bar{e}}} = M_{G}(X_{\bar{e}}=0)$.
\item
If $G=G_1 \vee G_2$, then we have
\[ \frac{M_G}{L_G} = \frac{M_{G_{1}}}{L_{G_{1}}} + \frac{M_{G_{2}}}{L_{G_{2}}}.
\]
\end{enumerate}
\end{lem}

\begin{Pf}
(1) is obvious from the definitions.
For (2), it is sufficient to show that
$L_{G_1} L_{G_2} = L_G$ and 
$M_{G_1} L_{G_2} + L_{G_1} M_{G_2} = M_G$,
but they are also obvious from the definitions and the above remarks.
\end{Pf}

Let $G$ be a hyperelliptic graph, and $D$ a polarization on $G$
with $\iota(D) = D$.
For any $\bar{e} \in {\mathrm{Ed}(G)}^{\sim}$,
$G^{\bar{e}}$ is simple and $D^{\bar{e}}= aP+bQ$,
where $P$ and $Q$ are the vertices and $a,b \in \mathbb{R}$.
Set $w(\bar{e})=\min\{a,b\}$.

The following theorem is a key result for our main theorem.

\begin{thm} \label{thm.b}
Let $G$ be a hyperelliptic graph, and $D$ a polarization given by
\[ D = \sum_{v \in \mathrm{Vert}_{\mathrm{n.f}}(G)}
(\nu_{G} (v) - 2) v +
\sum_{v' \in \mathrm{Vert}_{\mathrm{f}}(G)} a_{v'} v', 
\]
where $a_{v'} \in \mathbb{R}$.
Then, if $\mathrm{deg}(D)+2 \neq 0$, we have
\begin{align*} 
\epsilon(G,D) =
\sum_{\bar{e} \in {\mathrm{Ed}(G)}^{\sim}}
\Biggl(
\frac{2}{3} \frac{\mathrm{deg}(D)}{\mathrm{deg}(D)+2}
+ \frac{w(\bar{e})(\mathrm{deg}(D) - w(\bar{e}))}{\mathrm{deg}(D)+2}
\Biggr) X_{\bar{e}} 
+ \frac{2}{3} \frac{\mathrm{deg}(D)}{\mathrm{deg}(D)+2} \frac{M_G}{L_G},
\end{align*}
as rational functions on the length of each edge, i.e.,
on $\mathcal{M}({\mathrm{Ed}(G)}^{\sim})_{>0}$.
\end{thm}
In the rest of this section, we will give the proof of Theorem \ref{thm.b}.

\subsection{Preliminaries to the proof of Theorem \ref{thm.b}}

First of all, let us begin with direct calculations of
the admissible constants for 
$SG$ and $\mathbf{G}_{n-1}$ ($n>2$).

\begin{prop} \label{el.prop}
\begin{enumerate}
\renewcommand{\labelenumi}{(\arabic{enumi})}
\item 
Let $P$ and $Q$ be the two vertices of $SG$, and
$D= aP+bQ$ a polarization on $SG$.
Then, we have
\[ \epsilon(SG,D) =
\Biggl( \frac{2}{3} \frac{\mathrm{deg}(D)}{\mathrm{deg}(D) +2}
+ \frac{ab}{\mathrm{deg}(D) +2} \Biggr) X_{\bar{e}}. \]
\item
Let $Q$ be a non-fixed vertex of $\mathbf{G}_{n-1}$,
$\{ e_1, \ldots, e_n \}$ the set of edges which start from $Q$,
$P_i$ the other vertex of $e_i$ for $i=1,\ldots,n$, and
\[ D = (n-2)Q + (n-2) \iota(Q) + \sum_{i=1}^{n} a_{i} P_{i} \]
a hyperelliptic polarization on $\mathbf{G}_{n-1}$ with 
$\mathrm{deg}(D) +2 \neq 0$.
Then, we have
\begin{align*}
\epsilon (\mathbf{G}_{n-1},D) &=
\sum_{i=1}^{n} \Biggl(
\frac{2}{3} \frac{\mathrm{deg}(D)}{\mathrm{deg}(D) +2} +
\frac{a_i (\mathrm{deg}(D) - a_i)}{\mathrm{deg}(D) +2} \Biggr) X_{\bar{e}_{i}}
+ \frac{2}{3} \frac{\mathrm{deg}(D)}{\mathrm{deg}(D) +2}
\frac{(n-2) \sigma_n}{\sigma_{n-1}} ,
\end{align*}
where $\sigma_k$ is the $k$-th elementary symmetric polynomial on
$\{ X_{\bar{e}_{i}} \}$.
\end{enumerate}
\end{prop}

Note that $L_{\mathbf{G}_{n-1}} = \sigma_{n-1}$
and $M_{\mathbf{G}_{n-1}} = (n-2) \sigma_n$.

\begin{Pf}
For (1), it is easy to see by
\cite [Proposition 4.2, Corollary 4.3] {m-4}.
We will prove (2).

Let $\lambda = m_1 \bar{e}_{1}^{\ast}+ \cdots +m_n \bar{e}_{n}^{\ast}$ 
be a Lebesgue measure on
$\mathbf{G}_{n-1}$. Set
\begin{align*}
&\bar{\sigma}_k = \sigma_k (m_1,\ldots,m_n) \\
&\bar{\sigma}_{k}^{(i)} = \text{the $k$-th elementary symmetric polynomial}\\
&\quad \quad \quad \text{on
$\{ m_1,\ldots,m_{i-1},m_{i+1}, \ldots, m_n \}$.}
\end{align*}
By \cite [Lemma 3.7.] {z}, we have
\[ \mu_{(\mathbf{G}_{n-1};\lambda,D)}=
\frac{1}{\mathrm{deg}(D)+2} \Biggl(
\sum_{i=1}^{n} a_i \delta_i  +
\sum_{i=1}^{n} \frac{\bar{\sigma}_{n-2}^{(i)}}{\bar{\sigma}_{n_1}}
(de_i +d\iota(e_i)) \Biggr) . \]
Let 
\[ s_i : e_i \to [0,m_i] \]
be an arc-length parameter such that
$s_i (P_i) = 0$ and $s_i(Q) = m_i$.
We denote by the same symbol $s_i$ the parameter 
$s_i \circ \iota$ on $\iota(e_i)$.
Consider the following function on $\mathbf{G}_{n-1}$:
\[
g(x)=
\begin{cases}
\begin{split}
\displaystyle
&\frac{\bar{\sigma}_{n-2}^{(1)}}{2(\mathrm{deg}(D)+2)\bar{\sigma}_{n_1}}
s_1 (x)^2 \\
& \quad {} + \Biggl( \frac{a_1}{2(\mathrm{deg}(D)+2)} - 1 \Biggr) s_i (x) +
\gamma_1 
\end{split}
& \text{on $e_1$ or $\iota(e_1)$,} \\
\displaystyle
\begin{split}
& \frac{\bar{\sigma}_{n-2}^{(i)}}{2(\mathrm{deg}(D)+2)\bar{\sigma}_{n_1}}
s_i (x)^2 \\
& \quad {} + \frac{a_i}{2(\mathrm{deg}(D)+2)} s_i (x) +
\gamma_i
\end{split}
& \text{on $e_i$ or $\iota(e_i)$ for $i \neq 1$,}
\end{cases}
\]
where
\begin{align*}
\gamma_{1} &=
\frac{2}{3(\mathrm{deg}(D)+2)^2}
\Biggl( \bar{\sigma}_1 + (n-2)\frac{\bar{\sigma}_n}{\bar{\sigma}_{n-1}}
\Biggr) \\
& {} 
+ \frac{(\mathrm{deg}(D)-a_1)(\mathrm{deg}(D)-a_1+2)}{2(\mathrm{deg}(D)+2)^2}
m_1 
{} + \sum_{j > 1}
\frac{a_j (a_j + 2)}{2(\mathrm{deg}(D)+2)^2} m_j ,
\end{align*}
and
\begin{align*}
\gamma_i &=
\frac{2}{3(\mathrm{deg}(D)+2)^2}
\Biggl( \bar{\sigma}_1 + (n-2)\frac{\bar{\sigma}_n}{\bar{\sigma}_{n-1}}
\Biggr)
+ \sum_{j \neq 1,i}
\frac{a_j (a_j + 2)}{2(\mathrm{deg}(D)+2)^2} m_j \\
& {} +
\frac{a_1 - (a_1 + 1)(\mathrm{deg}(D) - a_1 + 2)}{(\mathrm{deg}(D)+2)^2} m_1 
+
\frac{a_i - (a_i + 1)(\mathrm{deg}(D) - a_i + 2)}{(\mathrm{deg}(D)+2)^2} m_i
\end{align*}
if $i \neq 1$.
Then, we can check by direct calculations that $g$ is continuous,
$\Delta(g) = \delta_{P_1} - \mu_{(G;\lambda,D)}$, and
$\int_G g \mu_{(G;\lambda,D)}=0$.
Thus, $g_{(G;\lambda,D)} (P_1,x) = g(x)$,
and by \cite [Lemma 4.1] {m-4},
we obtain the formula.
\end{Pf}

Following two lemmas are fundamental for hyperelliptic graphs.

\begin{lem} \label{mu-lem}
Let $(G,D)$ be a polarized hyperelliptic  graph with a Lebesgue measure 
$\lambda = \sum_{\bar{e}} l_{\bar{e}} \bar{e}^{\ast}$.
Then, the admissible metric is given by
\[ \mu_{(G;\lambda,D)} = \frac{1}{\mathrm{deg}(D)+2}
\Biggl( \delta_D - \delta_K 
+ \sum_{e \in \mathrm{Ed}(G)} 
\frac{P_{G}^{\bar{e}}(\lambda)}{L_{G}(\lambda)}de \Biggr) , \]
where $P_{G}^{\bar{e}}$ is the coefficient of $X_{\bar{e}}$ of $L_G$
when $L_G$ is regarded as a polynomial on $X_{\bar{e}}$.
\end{lem}

\begin{Pf}
We will prove the lemma by induction on the size of $G$.

If $\mathrm{sz}(G) \leq 2$, we obtain the formula by
direct calculations.

Now, suppose that $\mathrm{sz}(G) > 2$.
We may assume that $G$ is irreducible.
For any $\bar{e}_0 \in {\mathrm{Ed}(G)}^{\sim}$,
there exists a non-fixed vertex $v$ such that 
at least two one-jointed edges $e_1$ and $e_2$ start from $v$
and that $\bar{e}_1$ and $\bar{e}_2$ are different from $\bar{e}_0$.
Since this lemma is true for a hyperelliptic
graph $G_{\bar{e}}'$ if it is true for 
$G'$, we may assume that $\nu(v) = 3$,
i.e., exactly three edges start from $v$.
Note that
the third edge $\bar{e}_3$ which starts from $v$ is a disjoint edge
since $\mathrm{sz}(G) > 2$.
Let $H_1$ be the subgraph generated by
$\{ {e}_1,{e}_2, {e}_3 \} \cup 
\iota (\{ {e}_1,{e}_2, {e}_3 \})$ and
$G_1$ a metrized graph characterized by the following conditions.\\
(a) $\mathrm{Ed}(G_1) = \{ e_{1}',e_{1}'' \}$,
$e_{1}'$ and $e_{1}''$ are just the closed connected intervals,
and $e_{1}'$ intersects with $e_{1}''$ in one point $v_0'$.\\
(b) $l_{e_{1}'} = l_{e_{1}''} = 
\displaystyle l_{e_3}+\frac{l_{e_{1}}l_{e_{2}}}{l_{e_{1}}+l_{e_{2}}}$.\\
Let $v_{1}'$ (resp. $v_{1}''$) be the terminal point of 
$e_{1}'$ (resp. $e_{1}''$) which is not $v_0'$, and
$v_1$ the terminal point of $e_3$ which is not $v$.
We would like to consider another graph $G'$ which we obtain from $G$
by replacing the subgraph $H_1$ by $G_1$: the graph constructed in the
following way.\\
(1) Remove $H_1 \setminus \{ v_1, \iota (v_1) \}$ from $G$.\\
(2) Connect $G_1$ with $G \setminus (H_1 \setminus \{ v_1, \iota (v_1) \})$ 
by identifying $v_{1}'$ with
$v_1$, and $v_{1}''$ with $\iota(v_1)$. (See Figure \ref{change.of.graph}.)
\begin{figure}[htbp]
\begin{center}
\begin{picture}(200,100)
\put(20,20){\circle*{5}}
\put(10,50){\circle*{5}}
\put(30,50){\circle*{5}}
\put(20,80){\circle*{5}}
\put(60,20){\circle*{5}}
\put(60,50){\circle*{5}}
\put(60,80){\circle*{5}}
\put(60,20){\line(0,1){60}}
\put(20,20){\line(1,0){50}}
\put(20,80){\line(1,0){50}}
\put(10,50){\line(1,3){10}}
\put(10,50){\line(1,-3){10}}
\put(20,80){\line(1,-3){10}}
\put(20,20){\line(1,3){10}}
%\put(1,43){\small$v_{0}$}
\put(18,83){\small$v$}
\put(16,10){\small$\iota(v)$}
\put(50,84){\small$v_{1}$}
\put(50,10){\small$\iota(v_{1})$}
\put(5,65){\small$e_{1}$}
\put(25,65){\small$e_{2}$}
\put(34,83){\small$e_{3}$}
\put(90,47){$\longrightarrow$}
\put(140,50){\circle*{5}}
\put(170,20){\circle*{5}}
\put(170,80){\circle*{5}}
\put(170,50){\circle*{5}}
\put(140,50){\line(1,1){30}}
\put(140,50){\line(1,-1){30}}
\put(170,20){\line(0,1){60}}
\put(170,20){\line(1,0){10}}
\put(170,80){\line(1,0){10}}
\put(130,40){\small$v_{0}'$}
\put(145,68){\small$e_{1}'$}
\put(144,29){\small$e_{1}''$}
\put(153,84){\small$v_{1}' = v_{1}$}
\put(152,11){\small$v_{1}'' = \iota(v_{1})$}
\end{picture}
\caption{} \label{change.of.graph}
\end{center}
\end{figure}

$G'$ is again a hyperelliptic graph of size $(\mathrm{sz}(G)-1)$.
Note that we can naturally see each $e \in \mathrm{Ed}(G) \setminus 
\bigl( \{ e_1,e_2, e_3 \} \cup \iota (\{ e_1,e_2, e_3 \}) \bigr)$
as an edge of $G'$.
When we regard $e$ as an edge of $G'$, we denote it by $e'$.
Let $\lambda'$ be a Lebesgue measure on $G'$ such that
\begin{align*}
&l_{\bar{e}_{1}'} 
= l_{\bar{e}_3} 
+ \frac{l_{\bar{e}_1} l_{\bar{e}_2}}{l_{\bar{e}_1} + l_{\bar{e}_2}}, \\
&l_{\bar{e}'} = l_{\bar{e}} \quad \mbox{for } \bar{e}' \neq \bar{e}_{1}'.
\end{align*}
By the definition of $(G';\lambda')$, we have
\[ \frac{1}{l_{\bar{e}}+ r_{\bar{e}}}= 
\begin{cases}\displaystyle
\frac{1}{l_{\bar{e}'} + r_{\bar{e}'}} & 
\text{if $\bar{e} \neq \bar{e}_1, \bar{e}_2,\bar{e}_3$,}\\
\displaystyle
\frac{1}{l_{\bar{e}_1'} + r_{\bar{e}_1'}} &
\text{if $\bar{e} = \bar{e}_3$.} 
\end{cases}
\]
On the other hand,
we have
\[ \frac{2}{l_{\bar{e}_0'}+r_{\bar{e}_0'}}=
\frac{P_{G'}^{\bar{e}_0'}(\lambda')}{L_{G'}(\lambda')}, \]
by the induction hypothesis.
Therefore, the following suffices for our lemma
since 
$( X_{\bar{e}_1} + X_{\bar{e}_2} ) P_{G'}^{\bar{e}_0'} (\{ Y_{\bar{e}'} \}) 
= P_{G}^{\bar{e}_0}$ is automatic if the following is shown:
\[ ( X_{\bar{e}_1} + X_{\bar{e}_2} ) L_{G'} (\{ Y_{\bar{e}'} \}) = L_{G}, \]
where 
\[ Y_{\bar{e}'}=
\begin{cases}
X_{\bar{e}} & \text{if $\bar{e} \neq \bar{e}_1,\bar{e}_2,\bar{e}_3$}\\
X_{\bar{e}_3} + \displaystyle
\frac{X_{\bar{e}_1} X_{\bar{e}_2}}{X_{\bar{e}_1} + X_{\bar{e}_2}}
& \text{if $\bar{e}' = \bar{e}_1'$.}
\end{cases}
\]
This can be checked by direct calculations
if we use the second expression of $L_G$.
\end{Pf}

\begin{lem} \label{key-lem}
Let $G$ be a hyperelliptic graph and $D$ a polarization on $G$
such that the coefficient of
every non-fixed vertex $v$ in $D$ is $\nu_{G} (v) - 2$.
Then, there is a homogeneous polynomial $F$
of degree $(\mathrm{sz}(G) + 1)$ such that
\[ \epsilon(G,D) = \frac{F}{L_G}. \]
\end{lem}

\begin{Pf}
We may assume that $G$ is irreducible.
The following is a key claim.
\begin{cl} \label{cl.1}
In the same situation,
let $O$ be a fixed vertex, $e$ a one-jointed edge starting from $O$,
and $P_1$ the other terminal vertex of $e$.
Assume that there exists another one-jointed edge starting from $P_1$.
Then, for any vertex $P$, there exists a homogeneous polynomial $F_{P}$
of degree $(\mathrm{sz}(G) + 1)$ with
\[ g_{(G,D)}(O,P) = \frac{F_{P}}{L_G} . \]
\end{cl}
\begin{Pf}
We will show the claim by induction on $\mathrm{sz}(G)$.
For $\mathrm{sz}(G) = 2$, we have already obtained the claim in
Proposition \ref{el.prop}.
Assume that we have the claim for $\mathrm{sz}(G) \leq n$.
To simplify the notations, we only prove the claim for the graph like
f
Figure \ref{$G$}.
We can prove the claim for general hyperelliptic graphs in the same method.

\begin{figure}[htbp]
\begin{center}
\begin{picture}(260,100)
\put(20,20){\circle*{5}}
\put(10,50){\circle*{5}}
\put(30,50){\circle*{5}}
\put(20,80){\circle*{5}}
\put(50,20){\circle*{5}}
\put(50,50){\circle*{5}}
\put(50,80){\circle*{5}}
\put(80,20){\circle*{5}}
\put(80,50){\circle*{5}}
\put(80,80){\circle*{5}}
\put(140,20){\circle*{5}}
\put(140,50){\circle*{5}}
\put(140,80){\circle*{5}}
\put(170,20){\circle*{5}}
\put(160,50){\circle*{5}}
\put(180,50){\circle*{5}}
\put(170,80){\circle*{5}}
\put(50,20){\line(0,1){60}}
\put(80,20){\line(0,1){60}}
\put(140,20){\line(0,1){60}}
\put(20,20){\line(1,0){70}}
\put(20,80){\line(1,0){70}}
\put(10,50){\line(1,3){10}}
\put(10,50){\line(1,-3){10}}
\put(20,80){\line(1,-3){10}}
\put(20,20){\line(1,3){10}}
\put(125,80){\line(1,0){45}}
\put(125,20){\line(1,0){45}}
\put(160,50){\line(1,3){10}}
\put(160,50){\line(1,-3){10}}
\put(170,80){\line(1,-3){10}}
\put(170,20){\line(1,3){10}}
\put(100,47){$\cdots$}
\put(15,83){\small$P_{1}$}
\put(45,83){\small$P_{2}$}
\put(75,83){\small$P_{3}$}
\put(125,83){\small$P_{n-1}$}
\put(171,83){\small$P_{n}$}
\put(15,10){\small$\iota(P_{1})$}
\put(45,10){\small$\iota(P_{2})$}
\put(75,10){\small$\iota(P_{3})$}
\put(125,10){\small$\iota(P_{n-1})$}
\put(165,10){\small$\iota(P_{n})$}
\put(1,40){\small$O$}
\put(30,40){\small$Q_{1}$}
\put(52,40){\small$Q_{2}$}
\put(82,40){\small$Q_{3}$}
\put(115,40){\small$Q_{n-1}$}
\put(149,40){\small$Q_{n}$}
\put(179,40){\small$Q_{n+1}$}
\put(5,65){\small$e_{0}$}
\put(31,83){\small$e_{1}$}
\put(61,83){\small$e_{2}$}
\put(148,83){\small$e_{n-1}$}
\put(26,65){\small$f_{1}$}
\put(51,62){\small$f_{2}$}
\put(81,62){\small$f_{3}$}
\put(119,62){\small$f_{n-1}$}
\put(154,65){\small$f_{n}$}
\put(176,65){\small$f_{n+1}$}
\put(240,46){$\Bigg\updownarrow$}
\put(250,45){$\iota$}
\end{picture}
\caption{$G$} \label{$G$}
\end{center}
\end{figure}

In  virtue of \cite [Proposition 4.2] {m-4}, 
it is sufficient to show the case that
$D = P_1 + \cdots + P_n + \iota(P_1)  + \cdots + \iota(P_n)$.
Let us fix an arbitrary Lebesgue measure 
$\lambda$ on $G$ invariant under the involution, and fix
arc-length parameters $s_i$ and $t_j$ on $e_i$ and $f_j$ such that
$s_i (P_i)=0$, $s_i (P_{i+1})=l_i$, $t_j(Q_j)=0$, and $t_j(P_j)=m_j$,
where $l_i$ is the length of $e_i$ and $m_j$ is the length of $f_j$.
Let $\mu$ be the admissible metric of $(G;\lambda,D)$.
Set
\[g(x) = g_{(G;\lambda,D)}(O,x)=
\begin{cases}
\alpha_i s_i (x)^2 + \beta_i s_i (x) + \gamma_i & \mbox{on $e_i$,} \\
A_j t_j (x)^2 + C_j & \mbox{on $f_j$.}
\end{cases}
\]
We know by Lemma \ref{mu-lem} that $\alpha_i$ and $A_j$ is of form
$(\mathrm{poly.})/\bar{L}_G$, (i.e., for example,
there exists a homogeneous polynomial $H_{i}$ 
of degree $(\mathrm{sz}(G) - 1)$ determined by $(G,D)$ and $i$
such that $\alpha_i = H_{i}(\lambda)/L_G(\lambda)$.)
We can determine $\beta_i$'s inductively in the following way.
The first order differential equation at $P_n$ which comes from 
$\Delta g(x) = \delta_O - \mu$ gives
\[ -(2\alpha_{n-1} l_{n-1} + \beta_{n-1})
+(-2A_n m_n)+(-2A_{n+1} m_{n+1}) = 0, \]
hence we see $\beta_{n-1} = (\mathrm{poly.})/\bar{L}_G$.
Suppose that we know $\beta_k = (\mathrm{poly.})/ \bar{L}_G$ for $k \geq k_0$.
The first order differential equation at $P_{k_0}$ gives
\[ (-2\alpha_{k_0 -1} l_{k_0 -1} - \beta_{k_0 -1}) + (-2A_{k_0} m_{k_0})
+ \beta_{k_0} =0 , \]
hence we see $\beta_{k_{0} - 1} = (\mathrm{poly.})/\bar{L}_G$.
Thus, we have shown that all $\beta_i$'s are of form 
$(\mathrm{poly.})/\bar{L}_G$.

By the continuity of $g$, we have
\[
\begin{cases}
\alpha_{i-1} l_{i-1} ^2 + \beta_{i-1} l_{i-1} + \gamma_{i-1} = \gamma_{i}
& \mbox{for $i = 1, \ldots, n-1$,} \\
 A_j m_j ^2 + C_j = \gamma_j & \mbox{for $j = 1, \ldots, n-1$,} \\
A_{n} m_{n} ^2 + C_{n} = \alpha_{n-1} l_{n-1}^2 
+ \beta_{n-1} l_{n-1} + \gamma_{n-1} \\
A_{n+1} m_{n+1} ^2 + C_{n+1} =
\alpha_{n-1} l_{n-1}^2 + \beta_{n-1} l_{n-1} + \gamma_{n-1} ,
\end{cases}
\]
hence, in order to obtain the claim,
it is sufficient to show that $\gamma_0 = (\mathrm{poly.})/ \bar{L}_G$.

By the condition $\int_G g \mu =0$, we have
\begin{align*}
0 &=
2 \sum_{i=0}^{n-1} \Biggl(
\frac{1}{3} \alpha_{i} l_{i}^{2} + \frac{1}{2} \beta_{i} l_{i} + \gamma_{i}
\Biggr)(2\alpha_{i} l_{i})
+ 2 \sum_{j=1}^{n+1} \Biggl(
\frac{1}{3} A_{j} m_{j}^2 + C_{j} \Biggr)(2A_{j}m_{j}) \\
&= 2 \sum_{i=1}^{n-1} \Biggl(
\frac{1}{3} \alpha_{i} l_{i}^{2} + \frac{1}{2} \beta_{i} l_{i} + 
\sum_{k=1}^{i-1} (\alpha_{k} l_{k}^{2} + \beta_{k} l_{k})
+ \alpha_{0} l_{0}^{2} + \beta_{0} l_{0} +\gamma_{0} \Biggr)
(2 \alpha_{i} l_{i} ) \\
& \quad {}+ 2 \sum_{j=2}^{n} \Biggl(
\frac{1}{3} A_{j} m_{j}^2 + 
\sum_{k=1}^{j-1} (\alpha_{k} l_{k}^{2} + \beta_{k} l_{k}) - A_{j} m_{j}^{2}
+ \alpha_{0} l_{0}^{2} + \beta_{0} l_{0} +\gamma_{0} \Biggr)(2A_{j}m_{j}) \\
& \quad {}+ 2 \Biggl( \frac{1}{3} A_{n+1} m_{n+1}^2 +
\sum_{k=1}^{n-1} (\alpha_{k} l_{k}^{2} + \beta_{k} l_{k})
- A_{n+1} m_{n+1}^2 
+ \alpha_{0} l_{0}^{2} + \beta_{0} l_{0} +\gamma_{0} \Biggr)
(2A_{n+1} m_{n+1}) \\
& \quad {}+ 2 \Biggl(
\frac{1}{3} \alpha_{0} l_{0}^{2} + \frac{1}{2} \beta_{0} l_{0}
-(\alpha_{0} l_{0}^{2} + \beta_{0} l_{0})
+ \alpha_{0} l_{0}^{2} + \beta_{0} l_{0} + \gamma_{0} \Biggr)
(2\alpha_{0} l_{0}) \\
& \quad {}+ 2 \Biggl(
\frac{1}{3} A_{1} m_{1}^2 - A_{1} m_{1}^2 + 
\alpha_{0} l_{0}^{2} + \beta_{0} l_{0} + \gamma_{0} \Biggr)
(2A_{1} m_{1}) \\
&= (\alpha_{0} l_{0}^{2} + \beta_{0} l_{0} + \gamma_{0})
\Biggl( 2 \sum_{i=0}^{n-1} 2 \alpha_{i} l_{i}
+ 2 \sum_{j=1}^{n+1} 2A_{j}m_{j} \Biggr) \\
& \quad {}+ 2 \sum_{i=1}^{n-1} \Biggl(
\frac{1}{3} \alpha_{i} l_{i}^{2} + \frac{1}{2} \beta_{i} l_{i} + 
\sum_{k=1}^{i-1} (\alpha_{k} l_{k}^{2} + \beta_{k} l_{k})
\Biggr)(2\alpha_{i} l_{i}) \\
& \quad {}+ 2 \sum_{j=2}^{n} \Biggl(
- \frac{2}{3}A_{j} m_{j}^{2} + 
\sum_{k=1}^{j-1}(\alpha_{k} l_{k}^{2} + \beta_{k} l_{k}) \Biggr)
(2A_{j} m_{j}) \\
& \quad {}+ 2 \Biggl(
- \frac{2}{3} A_{n+1} m_{n+1}^2
+ \sum_{k=1}^{n-1}(\alpha_{k} l_{k}^{2} + \beta_{k} l_{k}) \Biggr)
(2A_{n+1} m_{n+1}) \\
& \quad {}+ 2 \Biggl(
- \frac{2}{3} \alpha_{0} l_{0}^{2} - \frac{1}{2}\beta_{0} l_{0} \Biggr)
(2\alpha_{0} l_{0})
+ 2\Biggl( - \frac{2}{3} A_{1} m_{1}^{2} \Biggr)(2A_{1} m_{1}) .
\end{align*}
By $\int_G \mu = 1$, we have
\[ 2 \Biggl( \sum_{i=0}^{n} 2 \alpha_{i} l_{i}
+ \sum_{j=1}^{n+1} 2A_{j}m_{j} \Biggr) = 1 \]
and thus,
\begin{align*}
-\gamma_{0} &=
2
\Biggl( - \frac{2}{3} \alpha_{0} l_{0}^{2} - \frac{1}{2}\beta_{0} l_{0} \Biggr)
(2\alpha_{0} l_{0})
+ 2 \Biggl( - \frac{2}{3} A_{1} m_{1}^{2} \Biggr)(2A_{1} m_{1}) \\
& \quad {} + 2 \sum_{i=1}^{n-1} \Biggl(
\frac{1}{3} \alpha_{i} l_{i}^{2} + \frac{1}{2} \beta_{i} l_{i} + 
\sum_{k=1}^{i-1} (\alpha_{k} l_{k}^{2} + \beta_{k} l_{k})
\Biggr)(2\alpha_{i} l_{i}) \\
& \quad {}+ 2 \sum_{j=2}^{n} \Biggl(
- \frac{2}{3}A_{j} m_{j}^{2} + 
\sum_{k=1}^{j-1}(\alpha_{k} l_{k}^{2} + \beta_{k} l_{k}) \Biggr)
(2A_{j} m_{j}) \\
& \quad {}+ 2 \Biggl(
- \frac{2}{3} A_{n+1} m_{n+1}^2
+ \sum_{k=1}^{n-1}(\alpha_{k} l_{k}^{2} + \beta_{k} l_{k}) \Biggr)
(2A_{n+1} m_{n+1})
+ (\alpha_{0} l_{0}^{2} + \beta_{0} l_{0}) .
\end{align*}
From the arguement so far,
we see that 
$g_{(G,D)}(O,O) 
= (\mbox{homog.poly.of deg $(2\mathrm{sz}(G) + 1)$.})/L_G ^2$ as rational
expressions on $\{ X_{\bar{e}} \}_{{\mathrm{Ed(G)}}^{\sim}}$.
Since $L_G$ is irreducible, it is enough for our claim to show that
$g_{(G,D)}(O,O) 
= (\mathrm{poly.})/(X_{\bar{e}_{0}}+ X_{\bar{f}_{1}})^a L_G$
for some nonnegative integer $a$.

Let $G'$ be the graph like Figure \ref{$G'$}.

\begin{figure}[htbp]
\begin{center}
\begin{picture}(260,100)
\put(20,20){\circle*{5}}
\put(10,50){\circle*{5}}
\put(20,80){\circle*{5}}
\put(50,20){\circle*{5}}
\put(50,50){\circle*{5}}
\put(50,80){\circle*{5}}
\put(80,20){\circle*{5}}
\put(80,50){\circle*{5}}
\put(80,80){\circle*{5}}
\put(140,20){\circle*{5}}
\put(140,50){\circle*{5}}
\put(140,80){\circle*{5}}
\put(170,20){\circle*{5}}
\put(160,50){\circle*{5}}
\put(180,50){\circle*{5}}
\put(170,80){\circle*{5}}
\put(50,20){\line(0,1){60}}
\put(80,20){\line(0,1){60}}
\put(140,20){\line(0,1){60}}
\put(20,20){\line(1,0){70}}
\put(20,80){\line(1,0){70}}
\put(10,50){\line(1,3){10}}
\put(10,50){\line(1,-3){10}}
\put(125,80){\line(1,0){45}}
\put(125,20){\line(1,0){45}}
\put(160,50){\line(1,3){10}}
\put(160,50){\line(1,-3){10}}
\put(170,80){\line(1,-3){10}}
\put(170,20){\line(1,3){10}}
\put(100,47){$\cdots$}
\put(15,83){\small$P_{1}$}
\put(45,83){\small$P_{2}$}
\put(75,83){\small$P_{3}$}
\put(125,83){\small$P_{n-1}$}
\put(171,83){\small$P_{n}$}
\put(15,10){\small$\iota(P_{1})$}
\put(45,10){\small$\iota(P_{2})$}
\put(75,10){\small$\iota(P_{3})$}
\put(125,10){\small$\iota(P_{n-1})$}
\put(165,10){\small$\iota(P_{n})$}
\put(1,38){\small$O'$}
\put(52,40){\small$Q_{2}$}
\put(82,40){\small$Q_{3}$}
\put(115,40){\small$Q_{n-1}$}
\put(149,40){\small$Q_{n}$}
\put(179,40){\small$Q_{n+1}$}
\put(5,65){\small$e_{0}'$}
\put(31,83){\small$e_{1}$}
\put(61,83){\small$e_{2}$}
\put(148,83){\small$e_{n-1}$}
\put(51,62){\small$f_{2}$}
\put(81,62){\small$f_{3}$}
\put(119,62){\small$f_{n-1}$}
\put(154,65){\small$f_{n}$}
\put(176,65){\small$f_{n+1}$}
\put(240,46){$\Bigg\updownarrow$}
\put(250,45){$\iota$}
\end{picture}
\caption{$G'$} \label{$G'$}
\end{center}
\end{figure}
Set
\[ D'= P_{1} + \cdots + P_{n}  + \iota(P_{1}) + \cdots + \iota(P_{n})+ 2 O' .\]
Let $\lambda'$ be the invariant measure 
by which the length of $e_{i}$ is $l_{i}$
for $i = 1, \ldots, n-1$, the length of $f_{j}$ is $m_{j}$
for $j = 2, \ldots, n+1$, and the length of $e_{0}'$ is 
$l_{0} m_{1}/(l_{0} + m_{1})$.
Set as before, for $i = 0, \ldots. n-1$ and ${j= 2, \ldots, n+1}$,
\[
g'(x) := g_{(G';\lambda',D')}(O',x) =
\begin{cases}
\alpha_{0}' {s_{0}'}(x)^{2} + \beta_{0}' s_{0}'(x) + \gamma_{0}' 
& \mbox{on $e_{0}'$,} \\
\alpha_{i}' {s_{i}'}(x)^{2} + \beta_{i}' s_{i}'(x) + \gamma_{i}' 
& \mbox{on $e_{i}$,} \\
A_{j}' {t_{j}'}(x)^{2} + C_{j}' & \mbox{on $f_{j}$.}
\end{cases}
\]
Of course, we have
$\alpha_{i}'=\alpha_{i}$ for $i \neq 0$ and $A_{j}'=A_{j}$
for $j \neq 1$,
and by the procedure in determining the $\beta_{i}$'s of $G$,
we also have $\beta_{i}'=\beta_{i}$ for $i \neq 0$.
Here, by the condition $\int_{G'} g \mu' = 0$, we have
\begin{align*}
0 &= 2 \sum_{i=1}^{n-1} \Biggl(
\frac{1}{3}\alpha_{i}' {l_{i}'}^{2} + \frac{1}{2} \beta_{i}' l_{i}' +
\gamma_{i}' \Biggr) (2 \alpha_{i}' l_{i}')
+ 2 \sum_{j=2}^{n+1} \Biggl( \frac{1}{3}A_{j}' {m_{j}'}^{2} 
+ C_{j}' \Biggr) (2 A_{j}'m_{j}') + \frac{2}{\deg D + 2} \gamma_{0}'\\
&= (\alpha_{0}' {l_{0}'}^{2} + \beta_{0}' l_{0}' + \gamma_{0}')
\Bigg( 2 \sum_{i=0}^{n-1} 2 \alpha_{i}' l_{i}' + 
2 \sum_{j=1}^{n+1} 2 A_{j}' m_{j}' + \frac{2}{\deg D +2} \Biggr) \\
& \quad {} + 2 \sum_{i=1}^{n-1} \Biggl(
\frac{1}{3} \alpha_{i} l_{i}^{2} + \frac{1}{2} \beta_{i} l_{i} + 
\sum_{k=1}^{i-1} (\alpha_{k} l_{k}^{2} + \beta_{k} l_{k})
\Biggr)(2\alpha_{i} l_{i}) \\
& \quad {}+ 2 \sum_{j=2}^{n} \Biggl(
- \frac{2}{3}A_{j} m_{j}^{2} + 
\sum_{k=1}^{j-1}(\alpha_{k} l_{k}^{2} + \beta_{k} l_{k}) \Biggr)
(2A_{j} m_{j}) \\
& \quad {}+ 2 \Biggl(
- \frac{2}{3} A_{n+1} m_{n+1}^2
+ \sum_{k=1}^{n-1}(\alpha_{k} l_{k}^{2} + \beta_{k} l_{k}) \Biggr)
(2A_{n+1} m_{n+1}) \\
& \quad {} - 2 \Biggl( \frac{2}{3} \alpha_{0}' {l_{0}'}^{2}
+ \frac{1}{2} \beta_{0}' l_{0}' \Biggr) (2 \alpha_{0}' l_{0}')
- \frac{2}{\deg D +2} (\alpha_{0}' {l_{0}'}^{2} + \beta_{0}' l_{0}'),
\end{align*}
hence, noting $\int_{G'} \mu' =1$ as before, we have
\begin{align*}
&2 \sum_{i=1}^{n-1} \Biggl(
\frac{1}{3} \alpha_{i} l_{i}^{2} + \frac{1}{2} \beta_{i} l_{i} + 
\sum_{k=1}^{i-1} (\alpha_{k} l_{k}^{2} + \beta_{k} l_{k})
\Biggr)(2\alpha_{i} l_{i}) \\
& \quad {}+ 2 \sum_{j=2}^{n} \Biggl(
- \frac{2}{3}A_{j} m_{j}^{2} + 
\sum_{k=1}^{j-1}(\alpha_{k} l_{k}^{2} + \beta_{k} l_{k}) \Biggr)
(2A_{j} m_{j}) \\
& \quad {}+ 2 \Biggl(
- \frac{2}{3} A_{n+1} m_{n+1}^2
+ \sum_{k=1}^{n-1}(\alpha_{k} l_{k}^{2} + \beta_{k} l_{k}) \Biggr)
(2A_{n+1} m_{n+1}) \\
&= 2 \Biggl( \frac{2}{3} \alpha_{0}' {l_{0}'}^{2}
+ \frac{1}{2} \beta_{0}' l_{0}' \Biggr) (2 \alpha_{0}' l_{0}')
+ \frac{2}{\deg D +2} (\alpha_{0}' {l_{0}'}^{2} + \beta_{0}' l_{0}')
-(\alpha_{0}' {l_{0}'}^{2} + \beta_{0}' l_{0}' + \gamma_{0}').
\end{align*}
Therefore,
\begin{align*}
- \gamma_{0} &=
2 
\Biggl( - \frac{2}{3} \alpha_{0} l_{0}^{2} - \frac{1}{2}\beta_{0} l_{0} \Biggr)
(2\alpha_{0} l_{0})
+ 2 \Biggl( - \frac{2}{3} A_{1} m_{1}^{2} \Biggr)(2A_{1} m_{1})
+ 2 \Biggl( \frac{2}{3} \alpha_{0}' {l_{0}'}^{2} 
+ \frac{1}{2} \beta_{0}' l_{0}' \Biggr) (2 \alpha_{0}' l_{0}') \\
& \quad {}+ \frac{2}{\deg D +2} 
(\alpha_{0}' {l_{0}'}^{2} + \beta_{0}' l_{0}')
-(\alpha_{0}' {l_{0}'}^{2} + \beta_{0}' l_{0}' + \gamma_{0}')
+ (\alpha_{0} {l_{0}}^{2} + \beta_{0} l_{0}) \\
&= - \frac{8}{3} (
\alpha_{0}^{2} l_{0}^{3} + A_{1}^{2} m_{1}^{3}
- {\alpha_{0}'}^{2} {l_{0}'}^{3} ) \\
& \quad {}
- \beta_{0} l_{0} (2\alpha_{0} l_{0})
+ \beta_{0}' l_{0}' (2\alpha_{0}' l_{0}')
+ \frac{2}{\deg D +2} (\alpha_{0}' {l_{0}'}^{2} + \beta_{0}' l_{0}') \\
& \quad {} -(\alpha_{0}' {l_{0}'}^{2} + \beta_{0}' l_{0}')
+ (\alpha_{0} {l_{0}}^{2} + \beta_{0} l_{0}) - \gamma_{0}'.
\end{align*}
$\beta_{0}$ (resp. $\beta_{0} '$) can be calculated with the first order 
differential equation at $O$ (resp. $O'$),
and this calculation shows that 
$\beta_{0}$ and $\beta_{0} '$ are just rational numbers
independent of the measure $\lambda$.
Noting that
\begin{align*}
\alpha_{0}' &= \frac{1}{(\deg(D) + 2)} \frac{1}{l_{e_{0}'} + r_{e_{0}'}} \\
&= \frac{1}{(\deg(D) + 2)} \frac{1}{l_{1} + r_{1}} = \alpha_{1},
\end{align*}
we see that 
\[ \alpha_{0}' = 
\frac{P_{G}^{\bar{f}_{1}}(\lambda)}{2(\deg(D) + 2)L_G (\lambda)}. \]
Now, look at the graph $G'$.
If we regard $e_{0}'' = e_{0}' \cup e_{1}$ 
and $\iota (e_{0}'') = \iota(e_{0}') \cup \iota(e_{1})$
as one edge, $G'$ is a hyperelliptic graph of size $\mathrm{sz}(G) - 1$.
Hence by the induction hypothesis, 
there exists a homogeneous polynomial $F_{1}$ on 
$\{ X_{\bar{e}} \}_{\bar{e} \in {\mathrm{Ed}(G')}^{\sim}}$
independent of $\lambda$ such that
$\gamma_{0}' = {F_{1}(\lambda')}/{L_{G'}(\lambda')}$.
As we saw in the proof of Lemma \ref{mu-lem},
\[ X_{\bar{e}_{0}''}= 
\frac{X_{\bar{e}_{0}} X_{\bar{f}_{1}}}{X_{\bar{e}_{0}} + X_{\bar{f}_{1}}}
+ X_{\bar{e}_{1}} \]
and
\[ (X_{\bar{e}_{0}} + X_{\bar{f}_{1}})L_{G'} = L_G,\] 
hence we have
a homogeneous polynomial $F_{2}$ on 
$\{ X_{\bar{e}} \}_{\bar{e} \in {\mathrm{Ed}(G)}^{\sim}}$
independent of $\lambda$ and a nonnegative integer $a$ independent of $\lambda$
such that
$\gamma_{0}' = {F_{2}(\lambda)}/{(l_{0} + m_{1})^{a} L_{G}(\lambda)}$.
Consequently, we see that all terms in
\begin{align*}
&- \frac{8}{3} (
\alpha_{0}^{2} l_{0}^{3} + A_{1}^{2} m_{1}^{3}
- {\alpha_{0}'}^{2} {l_{0}'}^{3} ) \\
&{} - \beta_{0} l_{0} (2\alpha_{0} l_{0})
+ \beta_{0}' l_{0}' (2\alpha_{0}' l_{0}')
+ \frac{2}{\deg D +2} (\alpha_{0}' {l_{0}'}^{2} + \beta_{0}' l_{0}') \\
&{} -(\alpha_{0}' {l_{0}'}^{2} + \beta_{0}' l_{0}')
+ (\alpha_{0} {l_{0}}^{2} + \beta_{0} l_{0}) - \gamma_{0}'
\end{align*}
but the first line are of 
form $(\mathrm{poly.})/(l_0 + m_1)^{a} \bar{L}_G$,
hence it suffices to show that the first line is also of that form,
i.e.,
there exist a homogeneous polynomial $F_{3}$ 
and a nonnegative integer $a$ such that
\[ \alpha_{0}^{2} l_{0}^{3} + A_{1}^{2} m_{1}^{3}
- {\alpha_{0}'}^{2} {l_{0}'}^{3} 
= \frac{F_{3}(\lambda)}{(l_{0} + m_{1})^{a} L_{G}(\lambda)} .
\]
Now for simplicity, set
$X_{0} = X_{\bar{e}_{0}}$, $X_{1} = X_{\bar{f}_{1}}$,
$P_{0} = P_{G}^{\bar{e}_{0}}$,
$P_{1} = P_{G}^{\bar{f}_{1}}$, and $P_{2} = P_{G}^{\bar{e}_{1}}$.
\[ \begin{cases} {\displaystyle
\sum_{i=0}^{n-1} 2 \alpha_{i} l_{i} +
\sum_{j=1}^{n+1} 2 A_{j} m_{j} = \frac{1}{2} }\\
{\displaystyle
\frac{1}{\deg D + 2} + \sum_{i=0}^{n-1} 2 \alpha_{i}' l_{i}' +
\sum_{j=2}^{n+1} 2 A_{j}' m_{j}'= \frac{1}{2}} ,
\end{cases} \]
which come from $\int_G \mu = 1$ and $\int_{G'} \mu ' = 1$,
give us
\[
2 \alpha_{0}' l_{0}' = -\frac{1}{\deg D + 2} + 2 A_{1} m_{1} + 
2 \alpha_{0} l_{0} . \]
Noting 
\begin{align*}
&\alpha_{0} = \frac{P_{0}(\lambda)}{2(\deg(D) + 2)L_{G}(\lambda)} \\
&A_{1} = \frac{P_{1}(\lambda)}{2(\deg(D) + 2)L_{G}(\lambda)} \\
&\alpha_{0}'= \alpha_{1} = \frac{P_{2}(\lambda)}{2(\deg(D) + 2)L_{G}(\lambda)},
\end{align*}
we have
\[ P_{2} X_{0} X_{1} \equiv
(X_{0} + X_{1})(P_{0} X_{0} + P_{1} X_{1}) \quad \mod L_G. \]
Thus, we are reduced to show
\[ (X_{0} + X_{1}) P_{0}^{2} X_{0}^{3}
+ (X_{0} + X_{1}) P_{1}^{2} X_{1}^{3}
- (P_{0} X_{0} + P_{1} X_{1})^{2} X_{0} X_{1}
\equiv 0 \quad \mod L_G \]
Let $B$ be a polynomial such that 
$X_{0} X_{1} B$ is the sum of all the monomials of $L_G$
which are divisible by $X_{0}X_{1}$.
Then, we have
\[ L_G = P_{0} X_{0} + P_{1} X_{1} - X_{0} X_{1} B \]
since any monomial in $L_G$ is divisible by $X_{0}$ or $X_{1}$.
Moreover, since
$ (P_{0} - X_{1} B)X_{0} = L_G - P_{1} X_{1}$, 
no monomials in $P_{0} - X_{1} B$ are divisible by $X_{0}$ or $X_{1}$.
Hence by symmetricity on $X_{0}$ and $X_{1}$,
$P_{1} - X_{0} B = P_{0} - X_{1} B$ and denote this
polynomial by $C$, which is free from $X_{0}$ and $X_{1}$.
Note that
\[ P_{0} X_{0} + P_{1} X_{1} \equiv X_{0} X_{1} B \quad \mod L_G \]
and
\begin{align*}
L_G &= P_{0} X_{0} + P_{1} X_{1} - X_{0} X_{1} B \\
&= X_{0} X_{1} B + X_{0} C + X_{0} X_{1} B + X_{1} C - X_{0} X_{1} B \\
&= X_{0} X_{1} B + (X_{0} +X_{1}) C .
\end{align*}
Then, we see
\begin{align*}
&(X_{0} + X_{1}) P_{0}^{2} X_{0}^{3}
+ (X_{0} + X_{1}) P_{1}^{2} X_{1}^{3}
- (P_{0} X_{0} + P_{1} X_{1})^{2} X_{0} X_{1}\\
= &(X_{0} + X_{1}) P_{0} X_{0}^{2} (P_{0} X_{0} + P_{1} X_{1})
+ (X_{0} + X_{1}) P_{1} X_{1} (P_{1} X_{1}^{2} - P_{0} X_{0}^{2})
-(P_{0} X_{0} + P_{1} X_{1})^{2} X_{0} X_{1}\\
\equiv &(X_{0} + X_{1}) P_{0} X_{0}^{2}(X_{0} X_{1} B) \\
&\quad{} +(X_{0} + X_{1}) P_{1} X_{1}
(X_{1}^{2} X_{0} B + X_{1}^{2} C - X_{0}^{2} X_{1} B - X_{0}^{2} C)
- (X_{0} X_{1} B)^{2} X_{0} X_{1} \quad \mod L_G \\
= &(X_{0} + X_{1}) P_{0} X_{0}^{2} (X_{0} X_{1} B)
+ (X_{1}^{2} - X_{0}^{2}) P_{1} X_{1} (X_{0} X_{1} B) \\
&\quad{}+ (X_{0} + X_{1}) (X_{1}^{2} - X_{0}^{2}) P_{1} X_{1} C
- (X_{0} X_{1} B)^{2} X_{0} X_{1}.
\end{align*}
Now
\begin{align*}
&(X_{0} + X_{1}) P_{0} X_{0}^{2} 
+ (X_{1}^{2} - X_{0}^{2}) P_{1} X_{1}
- (X_{0} X_{1} B) X_{0} X_{1} \\
= &(X_{0} + X_{1}) \bigl( X_{0} (X_{0} X_{1} B) + X_{0}^{2} C \bigr)
- (X_{0} X_{1} B) X_{0} X_{1} + (X_{1}^{2} - X_{0}^{2}) P_{1} X_{1} \\
= &X_{0}^{2} \bigl( X_{0} X_{1} B + (X_{0} +X_{1}) C \bigr)
+ (X_{1}^{2} - X_{0}^{2}) P_{1} X_{1} \\
\equiv &(X_{1}^{2} - X_{0}^{2}) P_{1} X_{1} \quad \mod L_G .
\end{align*}
Therefore, we have
\begin{align*}
&(X_{0} + X_{1}) P_{0}^{2} X_{0}^{3}
+ (X_{0} + X_{1}) P_{1}^{2} X_{1}^{3}
- (P_{0} X_{0} + P_{1} X_{1})^{2} X_{0} X_{1}\\
\equiv &(X_{1}^{2} - X_{0}^{2}) P_{1} X_{1} (X_{0} X_{1} B)
+ (X_{0} + X_{1}) (X_{1}^{2} - X_{0}^{2}) P_{1} X_{1} C
\quad \mod L_G \\
= &(X_{1}^{2} - X_{0}^{2}) P_{1} X_{1}
\bigl(X_{0} X_{1} B + (X_{0} +X_{1}) C \bigr) \\
\equiv &0 \quad \mod L_G,
\end{align*}
thus, we complete the proof of Claim \ref{cl.1}.
\end{Pf}

\begin{cl} \label{cl.2}
Let $O$ be a fixed vertex, $e$ an edge starting from $O$, and $P$ the other
terminal vertex of $e$.
Then, we have
\[ r_G (O,P) = \frac{\mbox{homog. poly. of deg $\mathrm{sz}(G)+1$}}{L_G}.  \]
\end{cl}
\begin{Pf}
Since $2/(X_{\bar{e}} +R_{e}) =P^{\bar{e}}_G / L_G$, we have
\[ R_{e} = \frac{2L_G}{P^{\bar{e}}_G} - X_{\bar{e}}.\]
Therefore,
\begin{align*} 
r_G (O,P) &= \frac{X_{\bar{e}} R_{e}}{X_{\bar{e}} + R_{e}} \\
&= \frac{(2L_G / P^{\bar{e}}_G)X_{\bar{e}} 
   - X_{\bar{e}} ^2}{2L_G / P^{\bar{e}}_G} \\
&= \frac{2L_G X_{\bar{e}} - X_{\bar{e}} ^2 P^{\bar{e}}_G}{2L_G},
\end{align*}
and thus, we obtain the claim.
\end{Pf}
\begin{cl} \label{cl.3}
Let $O$ be a fixed vertex as in Claim \ref{cl.1}. 
Then, for any vertex $P$, we have
\[ r_G (O,P) = \frac{\mbox{homog. poly. of deg $\mathrm{sz}(G)+1$}}{L_G}. \]
\end{cl}
\begin{Pf}
If $P$ is a fixed vertex, we can easily see that
$r_G (O,P)$ itself is a homogeneous polynomial of degree $1$.
In the argument below, hence,
we assume that $P$ is a non-fixed vertex.

We will prove the claim by induction on $\mathrm{sz}(G)$.

If $\mathrm{sz}(G)=2$, then we see 
that the claim is true by Claim \ref{cl.2} or
by direct calculations.

Suppose $\mathrm{sz}(G)>2$.
First we show the claim in the case where
there are at least two one-jointed edges
starting from $P$.
If $O'$ is the other terminal point of a one-jointed edge
starting from $P$, then 
\[ r_{G} (O',P) = \frac{\mbox{homog. poly. of deg $\mathrm{sz}(G)+1$}}{L_G} \]
by Claim \ref{cl.2}.
On the other hand, we know
\begin{align*}
& r_{G} (O,P) = g_{(G,D)} (O,O) - 2 g_{(G,D)} (O,P) + g_{(G,D)} (P,P) \\
& r_{G} (O',P) = g_{(G,D)} (O',O') - 2 g_{(G,D)} (O',P) + g_{(G,D)} (P,P), 
\end{align*}
hence we have
\[ r_{G} (O,P) -  r_{G} (O',P) = 
\frac{\mbox{homog. poly. of deg $\mathrm{sz}(G)+1$}}{L_G} \]
by Claim \ref{cl.1}.
Therefore, we obtain Claim \ref{cl.3}.
Next, we assume that there exists at most one one-jointed edge
starting from $P$.
Then, $G$ is not elementary and
we can find two distinct non-fixed vertices $P_{1}$ and $P_{2}$
whose classes in $\mathrm{Vert}(G)^{\sim}$ are
different from $P$, such that there exist at least two one-jointed
edges $e_{i,1}$ and $e_{i,2}$ starting from $P_{i}$ for $i = 1, 2$.
Let $O_{i}$ be the other terminal point of $e_{i,1}$ for $i = 1, 2$,
which is a non-fixed vertex.
Then, by the induction hypothesis and
the same argument in the proof of Lemma \ref{mu-lem}
or Claim \ref{cl.1}, we can see that there is a nonnegative integer $a'$ with
\[  r_{G} (O_{1},P) 
= \frac{\mbox{homog. poly. of deg 
$(\mathrm{sz}(G)+1+a')$}}{(X_{\bar{e}_{2,1}}+X_{\bar{e}_{2,2}})^{a'} L_G} . \]
Hence, again by Claim \ref{cl.1} and
\begin{align*}
& r_{G} (O,P) = g_{(G,D)} (O,O) - 2 g_{(G,D)} (O,P) + g_{(G,D)} (P,P) \\
& r_{G} (O_{1},P) 
= g_{(G,D)} (O_{1},O_{1}) - 2 g_{(G,D)} (O_{1},P) + g_{(G,D)} (P,P), 
\end{align*}
we see that there exist a nonnegative integer $a$
and a homogeneous polynomial $F_{1}$ 
of degree $(\mathrm{sz}(G)+1+a)$ which is coprime to 
$(X_{\bar{e}_{2,1}}+X_{\bar{e}_{2,2}})$, 
such that
\[ r_{G} (O,P)
= \frac{F_{1}}{(X_{\bar{e}_{2,1}}+X_{\bar{e}_{2,2}})^{a} L_G}. \]
In the same way, we see that
there exist a nonnegative integer $b$ and 
a homogeneous polynomial $F_{2}$ of degree $(\mathrm{sz}(G)+1+b)$ 
which is not divisible by 
$(X_{\bar{e}_{1,1}}+X_{\bar{e}_{1,2}})$,
such that
\[ r_{G} (O,P)
= \frac{F_{2}}{(X_{\bar{e}_{1,1}}+X_{\bar{e}_{1,2}})^{b} L_G} .  \]
Therefore,  we have
\[ (X_{\bar{e}_{1,1}}+X_{\bar{e}_{1,2}})^{b} F_{1}
= (X_{\bar{e}_{2,1}}+X_{\bar{e}_{2,2}})^{a} F_{2}, \]
and both $a$ and $b$ are equal to $0$ by the choice of $F_{1}$ and $F_{2}$.
\end{Pf}

In virtue of Claim \ref{cl.1}, Claim \ref{cl.3}
and \cite [Lemma 4.1] {m-4},
we obtain the Lemma \ref{key-lem}.
\end{Pf}

\subsection{Proof of Theorem \ref{thm.b}}
Now, we are ready to prove Thoerem \ref{thm.b}.
We will prove the theorem by induction on the size of $G$.

If $\mathrm{sz}(G)=1$ or $2$, it is Proposition \ref{el.prop}.
We assume $\mathrm{sz}(G)>2$.

\textit{Step 1.}
Suppose that $G$ is a one-point-sum of two graphs $G_1$ and $G_2$.
Let $D_i$ be the polarization 
$ D^{{\mathrm{Ed}(G_{i})}^{\sim}}$
on $G_i$ for $i=1,2$.
Then, $\mathrm{sz}(G_i) < \mathrm{sz}(G)$ for $i=1,2$,
 we obtain the formula
by the induction hypothesis and Lemma \ref{L-M-lem}.

\textit{Step 2.}
Suppose that $G$ is irreducible and $n= \mathrm{sz}(G)$.
Let us consider the following claim.
\begin{cl}
For $P \in S^{n+1}V$, we assume that (1) $P(X_{\bar{e}} =0) = 0$ for any 
$\bar{e} \in {\mathrm{Ed}_{1}(G)}^{\sim}$, and 
(2) $P(\{ X_{\bar{e}} = 0 \}_{\bar{e}  \in {\mathrm{Ed}_{0}(G)}^{\sim}} ) = 0$.
Then we have $P=0$.
\end{cl}

\begin{Pf}
Let $\lambda X_{\bar{e}_1} \cdots X_{\bar{e}_n} X_{\bar{e}_{n+1}}$
be a monomial in $P$.

If there is a disjoint edge class in $\{\bar{e}_1,\ldots, \bar{e}_{n+1} \}$,
then we have a one-jointed one which dose not appear in 
$\{\bar{e}_1,\ldots, \bar{e}_{n+1} \}$.
Therefore, we have $\lambda=0$ by  (1).

If there is no disjoint edge class in $\{\bar{e}_1,\ldots, \bar{e}_{n+1} \}$,
then all of them are one-jointed. Hence, we have also $\lambda=0$ by  (2).
\end{Pf}

Let us go back to the proof of Theorem \ref{thm.b}.
Set
\begin{align*} 
F =
& \sum_{\bar{e} \in {\mathrm{Ed}(G)}^{\sim}}
\Biggl(
\frac{2}{3} \frac{\mathrm{deg}(D)}{\mathrm{deg}(D)+2}
+ \frac{w(\bar{e})(\mathrm{deg}(D) - w(\bar{e}))}{\mathrm{deg}(D)+2}
\Biggr) X_{\bar{e}}
+ \frac{2}{3} \frac{\mathrm{deg}(D)}{\mathrm{deg}(D)+2} \frac{M_G}{L_G}.
\end{align*}
Since we know by Lemma \ref{key-lem}
\[ \epsilon(G,D) = \frac{P}{L_{G}} \]
for some homogeneous polynomial $P$ of degree $n+1$, in virtue of the claim,
it is sufficient to prove the following (1) and (2).
\begin{enumerate}
\renewcommand{\labelenumi}{(\arabic{enumi})}
\item
\begin{align*}
F (X_{\bar{e}_0}=0) = \epsilon(G_1 \vee G_2, D')
\end{align*}
for any $\bar{e}_0 \in {\mathrm{Ed}_{1}(G)}^{\sim}$,
where $G_1 \vee G_2 = G_{\bar{e}_0}$ and $D'=D_{\bar{e}_0}$.
\item
\begin{align*}
F (\{X_{\bar{e}}=0\}_{\bar{e} \in {\mathrm{Ed}_{0}(G)}^{\sim}}) 
= \epsilon(\mathbf{G}_{n}, D''),
\end{align*}
where $D''=(n-1)Q + (n-1) \iota(Q) + \sum_{i=1}^{n+1} a_{i} P_{i} $
in the same notation as Proposition \ref{el.prop}.
\end{enumerate}

(1) is equivalent to 
\[ \frac{M_G}{L_G} (X_{\bar{e}} = 0) = 
\frac{M_{G_1}}{L_{G_1}} + \frac{M_{G_2}}{L_{G_2}} \]
by \textit{Step 1},
which is nothing but Lemma \ref{L-M-lem}.
(2) is obvious from Proposition \ref{el.prop}.
Thus, we have achieved the conclusion.

\section{Proof of the main theorem}

In this section, we consider \emph{metrized} graphs only, hence
we denote the admissible constants by $\epsilon (G,D)$ instead of 
$\epsilon (G;\lambda,D)$.

We need one more lemma:
\begin{lem}
Let G be an irreducible hyperelliptic metrized graph with the Lebesgue measure
$\lambda = \sum_{\bar{e} \in 
{\mathrm{Ed}(G)}^{\sim}} l_{\bar{e}} \bar{e}^{\ast}$.
Then, we have
\[ \frac{\bar{M}_G}{\bar{L}_G} \leq
\sum_{\bar{e} \in {\mathrm{Ed}_{0}(G)}^{\sim}} l_{\bar{e}}
+ \frac{1}{4} \sum_{\bar{e} \in {\mathrm{Ed}_{1}(G)}^{\sim}} l_{\bar{e}}, \]
where $\bar{L}_G = L_G (\lambda)$ and $\bar{M}_G = M_G (\lambda)$.
Moreover, if $\mathrm{sz}(G) \leq 4$, then we have
\[ \frac{\bar{M}_G}{\bar{L}_G} \leq
\frac{1}{2} \sum_{\bar{e} \in {\mathrm{Ed}_{0}(G)}^{\sim}} l_{\bar{e}}
+ \frac{1}{4} \sum_{\bar{e} \in {\mathrm{Ed}_{1}(G)}^{\sim}} l_{\bar{e}}. \]
\end{lem}

\begin{Pf}
Let us consider the following:
\begin{align*}
\bar{L}_G \Biggl( \sum_{\bar{e} \in {\mathrm{Ed}_{0}(G)}^{\sim}} l_{\bar{e}}
+ \frac{1}{4} \sum_{\bar{e} \in {\mathrm{Ed}_{1}(G)}^{\sim}} l_{\bar{e}}
\Biggr)
= &\Biggl( \sum_{
\begin{smallmatrix}
\bar{e}_1, \ldots, \bar{e}_k \\
\mbox{\scriptsize {all distinct}}
\end{smallmatrix}}
\prod_{\bar{v}} \sigma_{(\bar{e}_1,\ldots,\bar{e}_k;\bar{v})}
l_{\bar{e}_1} \cdots l_{\bar{e}_k}
\Biggr)
\Biggl(\sum_{\bar{e} \in {\mathrm{Ed}_{0}(G)}^{\sim}} l_{\bar{e}} \Biggr) \\
&{}+
\Biggl( \sum_{
\begin{smallmatrix}
\bar{e}_1, \ldots, \bar{e}_k \\
\mbox{\scriptsize {all distinct}}
\end{smallmatrix}}
\prod_{\bar{v}} \sigma_{(\bar{e}_1,\ldots,\bar{e}_k;\bar{v})}
l_{\bar{e}_1} \cdots l_{\bar{e}_k}
\Biggr)
\Biggl(
\frac{1}{4} \sum_{\bar{e} \in {\mathrm{Ed}_{1}(G)}^{\sim}} l_{\bar{e}}
\Biggr),
\end{align*}
where $\bar{e}_1, \ldots, \bar{e}_k \in {\mathrm{Ed}_{0}(G)}^{\sim}$,
and $\bar{v}$ runs over all non-fixed vertices of 
$G_{{\mathrm{Ed}_{0}(G)}^{\sim} 
\setminus \{ \bar{e}_1, \ldots, \bar{e}_k \}}$.
Firstly, we would like to estimate the second line of the
right-hand-side of the above equality.
We kwow that an inequality
\begin{multline*}
(a_1 + \cdots + a_k)(a_1 a_2 \cdots a_{k-1} + a_2 a_3 \cdots a_{k} +
\cdots + a_{k} a_1 \cdots a_{k-2}) \geq
 k^2 a_1 a_2 \cdots a_{k} 
\end{multline*}
holds for positive numbers.
Hence, noting that ${\mathrm{Ed}_{1}(G)}^{\sim}
= {\mathrm{Ed}_{1}(G_{{\mathrm{Ed}_{0}(G)}^{\sim} \setminus
\{ \bar{e}_1,\ldots,\bar{e}_k \}})}^{\sim}$, we have
\begin{align*}
\Biggl( 
\prod_{\bar{v}} \sigma_{(\bar{e}_1,\ldots,\bar{e}_k;\bar{v})}
\Biggr)
\Biggl(
\frac{1}{4} \sum_{\bar{e} \in {\mathrm{Ed}_{1}(G)}^{\sim}} l_{\bar{e}}
\Biggr)
&=
\frac{1}{4}
\sum_{\bar{v}}
\Biggl(
\sigma_{(\bar{e}_1,\ldots,\bar{e}_k;\bar{v})}
\Biggl(
\sum_{
\begin{smallmatrix}
\bar{e} \in 
{\mathrm{Ed}_{1}(G)}^{\sim} \\
\mbox{\scriptsize {starts from $\bar{v}$}}
\end{smallmatrix}}
l_{\bar{e}}
\Biggr)
\prod_{\bar{v}' \neq \bar{v}} \sigma_{(\bar{e}_1,\ldots,\bar{e}_k;\bar{v}')}
\Biggr) \\
&\geq
\sum_{\bar{v}}
\Biggl(
\frac{\bigl(\nu^{1}(\bar{v})\bigr)^2}{4}
\tau_{(\bar{e}_1,\ldots,\bar{e}_k;\bar{v})}
\prod_{\bar{v}' \neq \bar{v}} \sigma_{(\bar{e}_1,\ldots,\bar{e}_k;\bar{v}')}
\Biggr),
\end{align*}
where $\nu^{1}(\bar{v}) = \nu_{(G)}^{1}(\bar{v})$
is the number of one-jointed edges of 
$G' = G_{{\mathrm{Ed}_{0}(G)}^{\sim} \setminus
\{ \bar{e}_1,\ldots,\bar{e}_k \}}$ starting from $v$.
We also write $\nu^{0}(\bar{v})$ and $\nu(\bar{v})$
omitting $G_{{\mathrm{Ed}_{0}(G)}^{\sim} \setminus
\{ \bar{e}_1,\ldots,\bar{e}_k \}}$.

Let us go on to the estimation of the first term.
For an arbitrary $\bar{e}_i \in \{ \bar{e}_1,\ldots,\bar{e}_k \}$,
let $v_{i,1}$ and $v_{i,2}$ be the terminal points of $e_i$,
and $v_{i}$ the vertex of $G_{{\mathrm{Ed}_{0}(G)}^{\sim} 
\setminus \{ \bar{e}_1, \ldots, \bar{e}_k \}}$ such that 
$v_{i,1}$ and $v_{i,2}$ go to $v_{i}$ when we contract $\bar{e}_{i}$.
Then, we can easily check that
\[
\sigma_{(\bar{e}_1,\ldots,\bar{e}_{i-1},\bar{e}_{i+1},\ldots,\bar{e}_k;
\bar{v}_i)}=
\tau_{(\bar{e}_1,\ldots,\bar{e}_k;\bar{v}_{i,1})}
\sigma_{(\bar{e}_1,\ldots,\bar{e}_k;\bar{v}_{i,2})}
+
\sigma_{(\bar{e}_1,\ldots,\bar{e}_k;\bar{v}_{i,1})}
\tau_{(\bar{e}_1,\ldots,\bar{e}_k;\bar{v}_{i,2})} .
\]
If we use this formula, we see that
\begin{align*}
&\sum_{i=1}^{k}
\Biggl(
\prod_{\bar{v}}
\sigma_{(\bar{e}_1,\ldots,\bar{e}_{i-1},\bar{e}_{i+1},\ldots,\bar{e}_k;
\bar{v})}
\Biggr) \\
=&
\sum_{i=1}^{k} \Biggl(
\tau_{(\bar{e}_1,\ldots,\bar{e}_k;\bar{v}_{i,1})}
\sigma_{(\bar{e}_1,\ldots,\bar{e}_k;\bar{v}_{i,2})}
\prod_{\bar{v} \neq \bar{v}_i}
\sigma_{(\bar{e}_1,\ldots,\bar{e}_{i-1},\bar{e}_{i+1},\ldots,\bar{e}_k;
\bar{v})} \\
&\quad \quad \quad \quad{} +
\sigma_{(\bar{e}_1,\ldots,\bar{e}_k;\bar{v}_{i,1})}
\tau_{(\bar{e}_1,\ldots,\bar{e}_k;\bar{v}_{i,2})}
\prod_{\bar{v} \neq \bar{v}_i}
\sigma_{(\bar{e}_1,\ldots,\bar{e}_{i-1},\bar{e}_{i+1},\ldots,\bar{e}_k;
\bar{v})}
\Biggr) \\
=&
\sum_{i=1}^{k} \Biggl(
\tau_{(\bar{e}_1,\ldots,\bar{e}_k;\bar{v}_{i,1})}
\sigma_{(\bar{e}_1,\ldots,\bar{e}_k;\bar{v}_{i,2})}
\prod_{\bar{v} \neq \bar{v}_{i,1},\bar{v}_{i,2}}
\sigma_{(\bar{e}_1,\ldots,\bar{e}_k;
\bar{v})} \\
&\quad \quad \quad \quad{} +
\sigma_{(\bar{e}_1,\ldots,\bar{e}_k;\bar{v}_{i,1})}
\tau_{(\bar{e}_1,\ldots,\bar{e}_k;\bar{v}_{i,2})}
\prod_{\bar{v} \neq \bar{v}_{i,1},\bar{v}_{i,2}}
\sigma_{(\bar{e}_1,\ldots,\bar{e}_k;
\bar{v})}
\Biggr) \\
=&
\sum_{i=1}^{k} \Biggl(
\tau_{(\bar{e}_1,\ldots,\bar{e}_k;\bar{v}_{i,1})}
\prod_{\bar{v} \neq \bar{v}_{i,1}}
\sigma_{(\bar{e}_1,\ldots,\bar{e}_k;
\bar{v})}
+
\tau_{(\bar{e}_1,\ldots,\bar{e}_k;\bar{v}_{i,2})}
\prod_{\bar{v} \neq \bar{v}_{i,2}}
\sigma_{(\bar{e}_1,\ldots,\bar{e}_k;
\bar{v})}
\Biggr)\\
=&
\sum_{\bar{v}} \Biggl(
\nu^{0} ( \bar{v})
\tau_{(\bar{e}_1,\ldots,\bar{e}_k;\bar{v})}
\prod_{\bar{v}' \neq \bar{v}}
\sigma_{(\bar{e}_1,\ldots,\bar{e}_k;
\bar{v}')}
\Biggr) .
\end{align*}
We have, therefore,
\begin{align*}
&\Biggl( \sum_{
\begin{smallmatrix}
\bar{e}_1, \ldots, \bar{e}_k \\
\mbox{\scriptsize {all distinct}}
\end{smallmatrix}}
\prod_{\bar{v}} \sigma_{(\bar{e}_1,\ldots,\bar{e}_k;\bar{v})}
l_{\bar{e}_1} \cdots l_{\bar{e}_k}
\Biggr)
\Biggl(\sum_{\bar{e} \in {\mathrm{Ed}_{0}(G)}^{\sim}} l_{\bar{e}} \Biggr) \\
=&
\sum_{\bar{e} \in {\mathrm{Ed}_{0}(G)}^{\sim}}
\Biggl(
\sum_{
\begin{smallmatrix}
\bar{e}_1, \ldots, \bar{e}_k \\
\mbox{\scriptsize {all distinct}}
\end{smallmatrix}}
\prod_{\bar{v}} \sigma_{(\bar{e}_1,\ldots,\bar{e}_k;\bar{v})}
l_{\bar{e}_1} \cdots l_{\bar{e}_k}
l_{\bar{e}}
\Biggr) \\
\geq&
\sum_{
\begin{smallmatrix}
\bar{e}_1, \ldots, \bar{e}_k \\
\mbox{\scriptsize {all distinct}}
\end{smallmatrix}}
\sum_{i=1}^{k}
\Biggl(
\prod_{\bar{v}}
\sigma_{(\bar{e}_1,\ldots,\bar{e}_{i-1},\bar{e}_{i+1},\ldots,\bar{e}_k;
\bar{v})}
\Biggr)
l_{\bar{e}_1} \cdots l_{\bar{e}_k}  \quad \quad (k \geq 1) \\
=&
\sum_{
\begin{smallmatrix}
\bar{e}_1, \ldots, \bar{e}_k \\
\mbox{\scriptsize {all distinct}}
\end{smallmatrix}}
\Biggl(
\sum_{j=1}^{k} \Biggl(
\nu^{0} ( \bar{v}_j )
\tau_{(\bar{e}_1,\ldots,\bar{e}_k;\bar{v}_{j})}
\prod_{\bar{v} \neq \bar{v}_{j}}
\sigma_{(\bar{e}_1,\ldots,\bar{e}_k;
\bar{v})}
\Biggr)
l_{\bar{e}_1} \cdots l_{\bar{e}_k}
\Biggr) .
\end{align*}
Thus, we have
\begin{align*}
&\bar{L}_G \Biggl(\sum_{\bar{e} \in {\mathrm{Ed}_{0}(G)}^{\sim}} l_{\bar{e}}
+ \frac{1}{4} \sum_{\bar{e} \in {\mathrm{Ed}_{1}(G)}^{\sim}} l_{\bar{e}}
\Biggr) \\
\geq&
\sum_{
\begin{smallmatrix}
\bar{e}_1, \ldots, \bar{e}_k \\
\mbox{\scriptsize {all distinct}}
\end{smallmatrix}}
\Biggl(
\sum_{\bar{v}} \Biggl(
\Biggl( \frac{\bigl(\nu^{1}(\bar{v})\bigr)^2}{4}+
\nu^{0} ( \bar{v} ) \Biggr)
\tau_{(\bar{e}_1,\ldots,\bar{e}_k;\bar{v})}
\prod_{\bar{v}' \neq \bar{v}}
\sigma_{(\bar{e}_1,\ldots,\bar{e}_k;
\bar{v}')}
\Biggr)
l_{\bar{e}_1} \cdots l_{\bar{e}_k}
\Biggr).
\end{align*}
Since
\[ \frac{\bigl(\nu^{1}(\bar{v})\bigr)^2}{4}+\nu^{0}( \bar{v})
\geq \nu^{1}(\bar{v}) + \nu^{0}( \bar{v}) - 2 
 = \nu( \bar{v}) - 2 \]
for any $v$, we obtain the first inequality.

We have also
\begin{align*}
&\bar{L}_G \Biggl( 
 \frac{1}{2} \sum_{\bar{e} \in {\mathrm{Ed}_{0}(G)}^{\sim}} l_{\bar{e}}
+ \frac{1}{4} \sum_{\bar{e} \in {\mathrm{Ed}_{1}(G)}^{\sim}} l_{\bar{e}}
\Biggr) \\
\geq&
\sum_{
\begin{smallmatrix}
\bar{e}_1, \ldots, \bar{e}_k \\
\mbox{\scriptsize {all distinct}}
\end{smallmatrix}}
\Biggl(
\sum_{\bar{v}} \Biggl(
\Biggl( \frac{\bigl(\nu^{1}(\bar{v})\bigr)^2}{4}+
\frac{\nu^{0}( \bar{v} ) }{2} \Biggr)
\tau_{(\bar{e}_1,\ldots,\bar{e}_k;\bar{v})}
\prod_{\bar{v}' \neq \bar{v}}
\sigma_{(\bar{e}_1,\ldots,\bar{e}_k;
\bar{v}')}
\Biggr)
l_{\bar{e}_1} \cdots l_{\bar{e}_k}
\Biggr).
\end{align*}
If $\mathrm{sz}(G) \leq 4$, then we see 
$\mbox{\#} {\mathrm{Ed}_{0}(G)}^{\sim} \leq 2$.
Therefore, we have
\[ \frac{\bigl(\nu^{1}(\bar{v})\bigr)^2}{4}+\frac{\nu^{0}( \bar{v})}{2}
\geq \nu^{1}(\bar{v}) + \nu^{0}( \bar{v}) - 2 
= \nu(\bar{v}) - 2 , \]
and we obtain the second inequality.
\end{Pf}

\begin{cor} \label{f-ineq}
In the same notation as that of \emph{Theorem \ref{thm.b}}, 
we have the following inequalities.
\begin{enumerate}
\renewcommand{\labelenumi}{(\arabic{enumi})}
\item
\begin{align*} 
\epsilon(G,D) \leq
\sum_{
w( \bar{e} ) \neq 0}
\Biggl(
\frac{4}{3} \frac{\mathrm{deg}(D)}{\mathrm{deg}(D)+2}
+ \frac{w(\bar{e})(\mathrm{deg}(D) - w(\bar{e}))}{\mathrm{deg}(D)+2}
\Biggr) l_{\bar{e}}
+ \sum_{
w( \bar{e} ) = 0} 
\frac{5}{6}  \frac{\mathrm{deg}(D)}{\mathrm{deg}(D)+2} l_{\bar{e}} 
\end{align*}
for any $G$ with a measure $\lambda = \sum l_{\bar{e}} \bar{e}^{\ast}$.
\item
\begin{align*} 
\epsilon(G,D) \leq
\sum_{
w( \bar{e} ) \neq 0}
\Biggl(
\frac{\mathrm{deg}(D)}{\mathrm{deg}(D)+2}
+ \frac{w(\bar{e})(\mathrm{deg}(D) - w(\bar{e}))}{\mathrm{deg}(D)+2}
\Biggr) l_{\bar{e}}
+ \sum_{
w( \bar{e} ) = 0} 
\frac{5}{6}  \frac{\mathrm{deg}(D)}{\mathrm{deg}(D)+2} l_{\bar{e}} 
\end{align*}
if every irreducible component of $G$ is of size less than $5$.
\end{enumerate}
\end{cor}

\begin{Pf}
Since we know that $e$ is one-jointed if $w(e)=0$,
this is immediate from Theorem \ref{thm.b} and the above lemma.
\end{Pf}

Let us start the proof of the main theorem. 
First of all, note the following fact
(cf. \cite [Theorem 5.6] {z}
\cite [Corollary 2.3] {m-1}
\cite [Theorem 2.1] {m-2}).
If $(\omega^{a}_{X/Y},\omega^{a}_{X/Y})_{a} > 0$, then we have
\[ \inf_{P \in \mathrm{Pic}^{0}(C)(\overline{K})} r_{C}(P)
\geq \sqrt{(g-1)(\omega^{a}_{X/Y},\omega^{a}_{X/Y})_{a}} ,\]
where $( , )_{a}$ is the admissible pairing.

Let $(G_y, \omega_y)$ be the 
polarized metrized graph by the configuration of $X_y$.
By the definition of the admissible pairing, we can see
\[ (\omega^{a}_{X/Y},\omega^{a}_{X/Y})_{a} =
(\omega_{X/Y},\omega_{X/Y}) - \sum_{y\in Y} \epsilon(G_y, \omega_y). \]

In virtue of \cite [Proposition 4.7] {c-h} and Noether's formula,
we have 
\begin{align*}
(\omega_{X/Y},\omega_{X/Y}) &=
\frac{g-1}{2g+1}\xi_{0}(X/Y) 
+ \sum_{j=1}^{[(g-1)/2]} 
\frac{6j(g-1-j)+2(g-1)}{2g+1}\xi_{j}(X/Y) \\
& \quad {} + \sum_{i=1}^{[ g/2 ]}
\Biggl( \frac{12i(g-i)}{2g+1} - 1 \Biggr) \delta_i (X/Y) .
\end{align*}
Let $(G_1,D_1)$ (resp. $(G_2,D_2)$)
be the polarized metrized graph obtained from $(G_y, \omega_y)$ by contracting 
all edges which correspond to nodes of positive type (resp. of type $0$).
Then, $G_1$ is a hyperelliptic graph 
as we saw in Example \ref{ex.of.h.g}
if we suitably redefine the set
of vertices and the set of edges, 
equipped with the involution invariant measure.
Moreover, the divisor $D_{1}$ is supported in the ``new''
set of vertices since the coefficient of each vertex corresponding to 
a $(-2)$-rational component is $0$.
On the other hand, $(G_2,D_2)$ is a tree.

Firstly, we will talk on the case of $g\geq 5$. 
By the definition of the polarized metrized dual graph,
Proposition \ref{sum-formula}
and Corollary \ref{f-ineq} (1), we have
\begin{align*} 
\epsilon(G_1,D_1)  &\leq
\frac{5(g-1)}{12g} \xi_{0}(X_{y})
+ \sum_{j=1}^{[(g-1)/2]}
\Biggl( 
\frac{4(g-1)}{3g}
+ \frac{2j(g-1-j)}{g}
\Biggr) \xi_{j}(X_{y}), \\
\epsilon(G_2,D_2)  &\leq
\sum_{i=1}^{[g/2]}
\Biggl( \frac{4i(g-1)}{g}-1 \Biggr)
\delta_i (X_{y}),
\end{align*}
and again by Proposition \ref{sum-formula}, we have
\begin{align*} 
\epsilon(G_y,\omega_y) \leq &
\frac{5(g-1)}{12g} \xi_{0}(X_{y})
+ \sum_{j=1}^{[(g-1)/2]}
\Biggl( 
\frac{4(g-1)}{3g}
+ \frac{2j(g-1-j)}{g}
\Biggr) \xi_{j}(X_{y})\\
&+ \sum_{i=1}^{[g/2]}
\Biggl( \frac{4i(g-1)}{g}-1 \Biggr)
\delta_i (X_{y}).
\end{align*}
Therefore, 
\begin{align*}
(\omega^{a}_{X/Y},\omega^{a}_{X/Y})_{a} & \geq
\frac{(g-1)(2g-5)}{12g(2g+1)} \xi_0 (X/Y) \\
& \quad {} + \sum_{j=1}^{[(g-1)/2]}
\frac{2(g-1)(3j(g-1-j)-g-2)}{3g(2g+1)} \xi_{j}(X/Y) \\
& \quad {} + \sum_{i=1}^{[g/2]}
\frac{4(g-1)i(g-i)}{g(2g+1)} \delta_{i}(X/Y).
\end{align*}
Now since $g \geq 5$, we have
\begin{align*}
3j(g-1-j)-g-2 &\geq 3(g-2)-g-2 \\
& = 2(g-4) \\
& > 0,
\end{align*}
which shows $(\omega^{a}_{X/Y},\omega^{a}_{X/Y})_{a} > 0$.
Thus, we obtain our theorem for $g \geq 5$.

Secondly, suppose  that $g \leq 4$.
If 
$\mathrm{sz}(G) >4$, 
then we can easily see that the degree of the polarization is larger than six.
On the other hand, the degree of the polarization is necessarily
equal to $2g-2 \leq 6$, which is a contradiction.
Therefore, we can use the inequality of  Corollary \ref{f-ineq} (2).
In the same way as above, we obtain an inequality

\begin{align*} 
\epsilon(G_y,\omega_y)  &\leq
\frac{5(g-1)}{12g} \xi_{0}(X_{y})
+ \sum_{j=1}^{[(g-1)/2]}
\Biggl( 
\frac{g-1}{g}
+ \frac{2j(g-1-j)}{g}
\Biggr) \xi_{j}(X_{y}) \\
& \quad {} +
\sum_{i=1}^{[g/2]}
\Biggl( \frac{4i(g-1)}{g}-1 \Biggr)
\delta_i (X_{y}),
\end{align*}
and 
\begin{align*}
(\omega^{a}_{X/Y},\omega^{a}_{X/Y})_{a} & \geq
\frac{(g-1)(2g-5)}{12g(2g+1)} \xi_0 (X/Y) \\
& \quad {} + \sum_{j=1}^{[(g-1)/2]}
\frac{(g-1)(2j(g-1-j)-1)}{g(2g+1)} \xi_{j}(X/Y) \\
& \quad {} + \sum_{i=1}^{[g/2]}
\frac{4(g-1)i(g-i)}{g(2g+1)} \delta_{i}(X/Y).
\end{align*}
Since all the coefficients of $\xi_{j}(X/Y)$ and $\delta_{i}(X/Y)$
are positive if $g \geq 3$,
we obtained the theorem for $g=3,4$.

\small{

}

\end{document}